\documentclass[reqno, oneside]{amsart}
\usepackage{style}

\begin{document}
\title[Dirac operators]{Szeg\H{o} condition  and scattering for one-dimensional Dirac operators}
\author{R.~V.~Bessonov}

\address{
\begin{flushleft}
bessonov@pdmi.ras.ru\\\vspace{0.1cm}
St.\,Petersburg State University\\  
Universitetskaya nab. 7/9, 199034 St.\,Petersburg, RUSSIA\\
\vspace{0.1cm}
St.\,Petersburg Department of Steklov Mathematical Institute of RAS\\
Fontanka 27, 191023 St.Petersburg,  RUSSIA\\
\end{flushleft}
}

\thanks{The author is supported by RFBR grant mol\_a\_dk 16-31-60053}

\begin{abstract}
We prove existence of modified wave operators for one-dimensional Dirac operators whose spectral measures belong to the Szeg\H{o} class on the real line.  
\end{abstract}

\maketitle

\section{Introduction}\label{s1}
\noindent Consider the one-dimensional Dirac operator $\Di_Q$ on the half-axis $\R_+ = [0, + \infty)$, 
\begin{equation}\label{do}
\Di_Q : X \mapsto JX' + Q X, \qquad J = \jm, \qquad Q = \left(\begin{smallmatrix}q_1 & q_2\\ q_2  & -q_1\end{smallmatrix}\right). 
\end{equation}
We assume that functions $q_1, q_2$ are real and absolutely integrable on compact subsets of $\R_+$. The ``free'' Dirac operator with potential $Q = 0$ will be denoted by $\Di_0$. Let $L^2(\R_+, \C^2)$ be the Hilbert space of measurable functions $X: \R_+ \to \C^2$ such that $\int_{\R_+}\|X(t)\|_{\C^2}^{2}\,dt < \infty$. The operator $\Di_Q$ defined by \eqref{do} on absolutely continuous functions $X \in L^2(\R_+, \C^2)$ such that $\Di_Q X \in L^2(\R_+, \C^2)$, $\langle X(0), \zo\rangle_{\C^2} = 0$, is the self-adjoint operator on $L^2(\R_+, \C^2)$, see Section 8.6 in \cite{LSb}. The standard wave operators for the pair $\Di_0$, $\Di_Q$ are defined as the limits
\begin{equation}\label{wo}
W_\pm(\Di_Q, \Di_0) =  \lim_{t \to \mp\infty} e^{it\Di_Q}e^{-it\Di_0} 
\end{equation}
in the strong operator topology in the case where these limits exist. Given a potential $Q = \left(\begin{smallmatrix}q_1 & q_2\\ q_2  & -q_1\end{smallmatrix}\right)$, we write $Q \in L^p$ if $q_1, q_2 \in L^p(\R_+)$. 
The classical scattering theory implies the existence of wave operators $W_\pm(\Di_Q, \Di_0)$ for  $Q \in L^1$, see, e.g., Theorem XI.9 in \cite{RSbook3}. In 2002, Christ and Kiselev \cite{ChK} proved the existence of wave operators $W_\pm(\Di_Q, \Di_0)$ for Dirac operators with potentials $Q \in L^p$, $1\le p < 2$. Using a different approach, Denisov \cite{Den02b} extended their result to the class $Q \in L^2$. His proof gives a formula for $W_\pm(\Di_Q, \Di_0)$ in terms of the Szeg\H{o} function of the spectral measure for $\Di_Q$. A measure $\mu = w\,dx + \mus$ on the real line $\R$ belongs to the Szeg\H{o} class $\sz$ if $(1+x^2)^{-1} \in L^1(\mu)$ and 
$$
\int_{\R}\frac{\log w(x)}{1+x^2}\,dx > -\infty,
$$
where $w$ is the density of the absolutely continuous part of $\mu$. The Szeg\H{o} function $D_\mu$ of $\mu \in \sz$ is the outer function in the upper half-plane $\C^+ = \{z \in \C: \; \Im z>0\}$ with modulus $\sqrt{w}$ on $\R$ and such that $D_\mu(i) >0$. In other words, 
\begin{equation}\label{eq94}
D_\mu(z) =  \exp\left(\frac{1}{\pi i}\int_{\R}\log \sqrt{w(x)}\left(\frac{1}{x-z} - \frac{x}{x^2+1}\right)dx\right), \qquad z \in \C^+.
\end{equation}
Let $\mu$ be the spectral measure of the operator $\Di_Q$, so that the generalized Fourier transform
\begin{equation}\label{ft}
\F_Q: X \mapsto \frac{1}{\sqrt{\pi}}\int_{\R_+}\langle X(t), \Psi(t, \bar z)\rangle_{\C^2}\,dt, \qquad z \in \C,
\end{equation} 
densely defined on functions with compact support, is the isometric operator from $L^2(\R_+, \C)$ to $L^2(\mu)$. Here 
$\left\langle \!\left(\begin{smallmatrix}a_1\\ a_2\end{smallmatrix}\right)\! , \!\left(\begin{smallmatrix}b_1\\ b_2\end{smallmatrix}\right) \!\right\rangle_{\C^2} = a_1 \bar b_1 + a_2 \bar b_2$, and $\Psi$ denotes the generalized eigenvector of $\Di_Q$: 
\begin{equation}\label{cp}
J\tfrac{\partial}{\partial t}\Psi(t, z) + Q \Psi(t, z) = z \Psi(t,z), \quad \Psi(0, z) = \oz, \qquad t \ge 0, \quad z \in \C.
\end{equation}
We choose the normalization in \eqref{ft} so that the Lebesgue measure on $\R$ is the spectral measure for $\Di_0$. It is known that every Dirac operator $\Di_Q$ on $\R_+$ with locally integrable potential $Q$ has the unique spectral measure. Moreover if $Q \in L^2$, then the spectral measure $\mu$ of $\Di_Q$ belongs to the Szeg\H{o} class $\sz$, see Theorem~12.1 and Theorem 14.4 in \cite{Den06}. It was proved by Denisov \cite{Den02b} that for all $Q \in L^2$ the wave operators $W_\pm(\Di_Q, \Di_0)$ exist and 
\begin{equation}\label{eq29}
W_-(\Di_Q, \Di_0) = \gamma\F_Q^{-1}\chi_{E}D_\mu^{-1}\F_0, \qquad W_+(\Di_Q, \Di_0) = \bar\gamma\F_Q^{-1}\chi_{E}\ov{D_\mu^{-1}}\F_0, 
\end{equation}
where $\chi_E$ is the indicator of a Borel set $E$ such that $|\R\setminus E| = 0$, $\mus(\R\setminus E) = 0$. Here and in what follows $\mus$ denotes the singular part of $\mu$. We also denoted by $\chi_{E}D_\mu^{-1}$ the multiplication operator from $L^2(\R)$ to $L^2(\mu)$ taking $f$ into $\chi_{E}D_\mu^{-1} f$. The parameter $\gamma \in \C$ in \eqref{eq29} depends only on $Q$, $|\gamma| = 1$. Moreover, $\gamma = 1$ for potentials $Q$ with zero entries on the diagonal ($q_1 = 0$). It is not difficult to see that the operator in \eqref{eq29} is unitary from $L^2(\R_+, \C^2)$ onto the absolutely continuous subspace of $\Di_Q$. In other words, the operators $W_\pm(\Di_Q, \Di_0)$ for $Q \in L^2$ are complete.
 
\medskip

It is known that for every $p>2$ there exists a potential $Q \in L^p$ such that the absolutely continuous spectrum of $\Di_Q$ is empty. In particular, the strong wave operators $W_\pm(\Di_Q, \Di_0)$ for such potentials $Q$ do not exist. Thus, the result of Denisov is sharp in the scale of $L^p$-spaces. However, there are locally integrable potentials $Q \notin \cup_{1 \le p \le 2} L^p$ such that the spectral measures of operators $\Di_Q$ belong to the Szeg\H{o} class $\sz$. For such measures $\mu$ the Szeg\H{o} function $D_\mu$ and the right hand side of \eqref{eq29} are well-defined. It is natural to expect that  the limits in \eqref{wo}  defining wave operators $W_\pm(\Di_Q, \Di_0)$ exist as well. Examples constructed by Teplyaev \cite{Tep05} show that this is not the case, see the discussion after Theorem 14.6 in \cite{Den06}. In fact, the number $\gamma$ in \eqref{eq29} is the limit of a function $e^{i\phi(t)}$ related to the solution of \eqref{cp} at the point $z = i$. If $Q \notin L^2$, this function can have non-vanishing oscillations as $t \to \infty$. To handle these oscillations, one can introduce the modified wave operators 
\begin{equation}\label{eq56}
W^m_\pm(\Di_Q, \Di_0) =  \lim_{t \to \mp\infty} e^{it\Di_Q}\Sigma_\phi e^{-it\Di_0}, \qquad \Sigma_\phi X = \rrot X. 
\end{equation}
Following Denisov's approach for $Q \in L^2$, we prove existence of $W^m_\pm(\Di_Q, \Di_0)$ for general Dirac operators with spectral measures in Szeg\H{o} class. This requires some new constructions and estimates inspired by the theory of orthogonal polynomials on the unit circle. Our main result can be summarized as follows.
\begin{Thm}\label{t1}
Let $q_1$, $q_2$ be real-valued functions on $\R_+$ such that $q_1, q_2 \in L^1[0,r]$ for every $r>0$, and let $Q = \left(\begin{smallmatrix}q_1 & q_2\\ q_2 & -q_1\end{smallmatrix}\right)$. Assume that the spectral measure $\mu$ of $\Di_Q$ belongs to $\sz$. Then there exists a function $\phi$ such that the strong wave operators $W^m_\pm(\Di_Q, \Di_0)$ in \eqref{eq56} exist and complete. Moreover, $W^m_\pm(\Di_Q, \Di_0)$ coincide with the operators in \eqref{eq29} for $\gamma  = 1$. If $q_1 = 0$, then $\phi = 0$ and the usual wave operators $W_\pm(\Di_Q, \Di_0)$ exist and complete.
\end{Thm}
Applying results of \cite{BD2017}, one can explicitly characterize potentials $Q = \left(\begin{smallmatrix}0 & q\\ q  & 0\end{smallmatrix}\right)$ such that the spectral measure of $\Di_Q$ belongs to the Szeg\H{o} class. For $n \ge 0$, introduce the functions
$$
g_n(t) = \exp\left(2\int_{n}^{t}q(s)\,ds\right), \qquad t \in [n, n+2). 
$$
As we will see in Section \ref{s4}, the spectral measure of $\Di_Q$ belongs to $\sz$ if and only if 
\begin{equation}\label{eq57}
\sum\limits_{n \ge 0} \left(\int_{n}^{n+2}g_n(t)\,dt\cdot\int_{n}^{n+2}\frac{dt}{g_n(t)} - 4\right) < \infty. 
\end{equation} 
Taylor expansion of the function $e^x$ at $x = 0$ shows that if 
\begin{equation}\label{eq54}
\sum_{n \ge 0} \max_{n \le t \le n+1}\left(\int_{n}^{t}q(s)\,ds\right)^2 < \infty,
\end{equation}
then relation \eqref{eq57} holds. If one replaces $q$ by $|q|$ in \eqref{eq54}, then the resulting class of potentials will coincide with $\ell^2(L^1(\R_+))$, the endpoint of the scale $\ell^p(L^1(\R_+))$, $1\le p<2$, treated by Christ and Kiselev in \cite{ChK}. It is interesting to note that the class of potentials described by \eqref{eq54} contains some non-decaying potentials. For example, the function $q: x \mapsto \sin(x^2)$ satisfies \eqref{eq54} but does not belong to any class $L^p(\R_+)$ or $\ell^p(L^1(\R_+))$, $1 \le p<\infty$. Theorem \ref{t1} implies the following result.
\begin{Cor}\label{c2} Let $q$ be a function on $\R_+$ such that $q \in L^1[0,r]$ for all $r>0$,
and let $\Di_Q$ be the Dirac operator on $L^2(\R_+, \C^2)$ with the potential $Q = \left(\begin{smallmatrix}0 & q\\ q  & 0\end{smallmatrix}\right)$. If $q$ satisfies \eqref{eq57} or \eqref{eq54}, then the wave operators $W^\pm(\Di_Q, \Di_0)$ exist. 
\end{Cor}
A part of our analysis is close to Khrushchev ideas in the theory of orthogonal polynomials on the unit circle. One of the central results of Khrushchev paper \cite{KH01} concerns a logarithmic convergence of orthogonal polynomials, see Theorem 2.5 in \cite{KH01}. 
Let us formulate its analogue for measures $\mu$ on the real line.   
\begin{Thm}\label{t2}
Let $\mu = w\,dx + \mus$ be a measure in the Szeg\H{o} class $\sz$, and let $\{\B_{r}\}$ be the chain de Branges spaces isometrically embedded into $L^2(\mu)$. Then  
$$
\lim_{r \to \infty}\int_{\R}\left|\log\frac{1}{|E_r(x)|^2} - \log w(x)\right|\frac{dx}{1+x^2} = 0,
$$
for some Hermite-Biehler functions $E_r$ generating $\B_r$.
\end{Thm}

More details on Theorem \ref{t2} can be found in Section \ref{s3}. Section \ref{s5} deals with the main instrument of the present paper -- regularized Krein's orthogonal entire functions. In Section \ref{s4} we prove Theorem \ref{t1}. Next section concerns the theory of canonical Hamiltonian systems that will be used in the paper. 

\section{Canonical Hamiltonian systems}\label{s2}
In this section we collect some known facts related to the spectral theory of canonical Hamiltonian systems, discuss the definition and properties of an entropy function introduced in \cite{BD2017}, recall the notion of de Branges chains, and reduce consideration of Dirac systems to canonical Hamiltonian systems. 

\subsection{Canonical Hamiltonian systems} The canonical Hamiltonian system on the positive half-axis $\R_+ = [0,+\infty)$ is the differential equation
\begin{equation}\label{cs}
\begin{cases}
J \tfrac{\partial}{\partial t} M(t,z) = z \Hh(t) M(t,z), \\ M(0, z) = \idm,
\end{cases} 
\quad t \in \R_+, \quad z \in \C.
\end{equation}
As before, $J = \jm$ is the constant sign matrix, the Hamiltonian $\Hh$ is a matrix-valued mapping such that
$$
\Hh = \sth, \quad \trace\Hh(t) >0, \quad \det \Hh(t) \ge 0, \quad t \in\R_+.
$$
The functions $h_1$, $h_2$, $h$ are assumed to be real-valued and belong to $L_{loc}^1(\R_+)$. Here and below $L_{loc}^1(\R_+)$ denotes the set of functions $f$ on $\R_+$ such that $f \in L^1[0,r]$ for every $r>0$. We will often represent the solution $M$ of \eqref{cs} in the form 
\begin{equation}\label{eq32}
M(t, z) = \begin{pmatrix}\Theta(t,z), \Phi(t,z)\end{pmatrix} = \mm, \qquad t \in \R_+, \quad z \in \C.
\end{equation} 
Since $\Hh$ is locally integrable, the function $M$ is locally absolutely continuous with respect to $t$ if the spectral parameter $z \in \C$ is fixed. It is also easy to see that for every $t\in \R_+$ the entries of $M$ are entire functions with respect to $z$. A Hamiltonian $\Hh$ on $\R_+$ is called singular if
\begin{equation}\label{eq33}
\int_{0}^{+\infty}\trace\Hh(t)\, dt = + \infty.
\end{equation}
We say that $\Hh$ is trivial if $\Hh$ coincides with one of two matrices 
$\left(
\begin{smallmatrix}
1 & 0\\
0 & 0\\
\end{smallmatrix}
\right)$, 
$\left(
\begin{smallmatrix}
0 & 0\\
0 & 1\\
\end{smallmatrix}
\right)$ almost everywhere on $\R_+$. Fix a parameter $\omega \in \R \cup \{\infty\}$. The Weyl function of $\Hh$ is defined by
\begin{equation}\label{mf}
m(z) = \lim_{t \to +\infty}\frac{\omega\Phi^+(t,z) + \Phi^-(t,z)}{\omega\Theta^+(t,z) + \Theta^-(t,z)}, \qquad z \in \C^+,
\end{equation}
where $\C^+ = \{z \in \C^+:\; \Im z >0\}$. It can be shown \cite{HSW} that for every singular nontrivial Hamiltonian $\Hh$ the expression under the limit in \eqref{mf} is correctly defined for large $t>0$ (the denominator is non-zero), it does not depend on $\omega$, and, moreover, $m$ is the analytic function in $\C^+$ with strictly positive imaginary part. In particular, $m$ admits the Herglotz representation 
\begin{equation}\label{hg}
m(z) = \frac{1}{\pi}\int_{\R}\left(\frac{1}{x - z} - \frac{x}{x^2+1}\right) d\mu(x) + bz + a, \quad z \in\C^+,
\end{equation}
where $\mu$ is a Radon measure on $\R$ such that $\int_{\R}\frac{d\mu(x)}{x^2+1} < \infty$, $b \ge 0$, and $a \in \R$. Krein -- de Branges theorem \cite{KK68}, \cite{dbbook} says that any function with positive imaginary part is the Weyl function of a singular nontrivial Hamiltonian on $\R_+$, see also \cite{Winkler95}, \cite{Romanov}. The measure $\mu$ in \eqref{hg} is called the spectral measure of the Hamiltonian $\Hh$. We will say that $\mu$ generates a singular nontrivial Hamiltonian $\Hh$ on $\R_+$ if the Weyl function $m$ of $\Hh$ satisfies \eqref{hg} for $\mu$ and $a=b=0$. A $\sigma$-finite Borel measure $\mu$ on $\R_+$ is called even if $\mu(I) = \mu(-I)$ for every interval $I$. It is well-known that even measures generate diagonal Hamiltonians, see, e.g., \cite{BD2017}. By the diagonal Hamiltonian $\diag (h_1, h_2)$ we mean the matrix-function $\left(\begin{smallmatrix}h_1&0\\0&h_2\end{smallmatrix}\right)$ on $\R_+$. 

\medskip

The main result of \cite{BD2017} is the following Szeg\H{o}-type theorem for diagonal canonical Hamiltonian systems. 
\begin{Thm}\label{t4} Let $\mu$ be an even measure on $\R_+$ such that $(1+x)^{-2} \in L^1(\mu)$, and let 
$\Hh = \diag(h_1, h_2)$ be a Hamiltonian generated by $\mu$. Then $\mu \in \sz$ if and only if $\sqrt{\det \Hh} \notin L^1(\R_+)$ and 
\begin{equation}\label{eq34}
\sum_{n = 0}^{+\infty}  \left(\int_{\eta_n}^{\eta_{n+2}}h_1(s)\,ds \cdot \int_{\eta_n}^{\eta_{n+2}}h_2(s)\,ds - 4\right) < \infty,
\end{equation}
where $\eta_n = \min\left\{t \ge 0: \int_{0}^{t}\sqrt{\det\Hh(s)}\,ds = n\right\}$ for $n \ge 0$. 
\end{Thm}
As we will see in Section \ref{s24}, Theorem \ref{t4} implies that the spectral measure of the Dirac operator $\Di_Q$ with a potential $Q = \left(\begin{smallmatrix}0 & q\\ q  & 0\end{smallmatrix}\right)$ belong to the Szeg\H{o} class $\sz$ if and only if relation \eqref{eq57} holds. 

\medskip
 
\subsection{Entropy function of a Hamiltonian}\label{s22} Let $\Hh$ be a singular nontrivial Hamiltonian on $\R_+$. For every $r \ge 0$ define $\Hh_r$ to be the Hamiltonian on $\R_+$ taking $x$ into $\Hh(x+r)$. Let $m_r$, $\mu_r$, $b_r$, $a_r$ denote the Weyl function, the spectral measure, and the coefficients in Herglotz representation \eqref{hg} for $m_r$. Define 
\begin{align*}\label{eq63}
\I_{\Hh}(r) &= \Im m_r(i) = \frac{1}{\pi}\int_{\R}\frac{d\mu_r(x)}{1+x^2} + b_r, \\
\Rr_{\Hh}(r) &= \Re m_r(i) = a_r,\\
\J_{\Hh}(r) &= \frac{1}{\pi}\int_{\R}\frac{\log w_r(x)}{1+x^2}\,dx, 
\end{align*}  
where $w_r$ is the density of the absolutely continuous part of $\mu_r = w_r \,dx + \mu_{r,\mathbf{s}}$. It can be shown that $\Rr_{\Hh}$ is identically zero if the Hamiltonian $\Hh$ is diagonal, see Lemma 2.2 in \cite{BD2017}. Following \cite{BD2017}, define the entropy function of $\Hh$ by
$$
\K_\Hh(r) = \log \I_{\Hh}(r) - \J_{\Hh}(r), \qquad r \ge 0. 
$$
Since  $b_r \ge 0$ by construction, from Jensen inequality we see that $\K_\Hh(r)\ge 0$. 
Consider the Hamiltonian 
\begin{equation}\label{cook59}
\widehat\Hh_r(t) =
\begin{cases}
\Hh(t), \quad &t \in [0,r), \\
\left(\begin{smallmatrix}
c_1(r) & c(r)\\
c(r) & c_2(r)
\end{smallmatrix}\right), \quad  &t \in [r, + \infty),
\end{cases}
\end{equation}
where $c_1(r) = 1/\I_{\Hh}(r)$, $c(r) = \Rr_{\Hh}(r)/\I_{\Hh}(r)$,  $c_2(r) = (\I_{\Hh}^2(r) + \Rr_{\Hh}^2(r))/\I_{\Hh}(r)$.
The Hamiltonian $\widehat\Hh_r$ coincides with $\Hh$ on $[0,r)$, is constant on $[r,+\infty)$, and, moreover, we have $\K_{\widehat\Hh_r}(0) \le \K_{\Hh}(0)$, see Appendix. We call $\widehat\Hh_r$ the Bernstein-Szeg\H{o} approximation to $\Hh$. The following three results were proved in \cite{BD2017} for diagonal Hamiltonians, see Appendix for the general case.  
\begin{Lem}\label{l1}
Let $\Hh$ be a singular nontrivial Hamiltonian on $\R_+$ and let $\mu = w\,dx + \mus$ be the spectral measure of $\Hh$. 
Assume that $\mu \in \sz$. Then for every $r \ge 0$ the measure $\mu_r = w_r\,dx + \mu_{r, \mathbf{s}}$ is in $\sz$ and we have 
\begin{itemize}
\item[$(a)$] $m(z) = \frac{G_r(z)}{F_r(z)}$ for all $z \in \C^+$, 
\item[$(b)$] $w(x) = \frac{w_r(x)}{|F_r(x)|^{2}}$ for almost all $x \in \R$,
\end{itemize}
where $G_r: z \mapsto \Phi^+(r,z) + m_r(z)\Phi^-(r,z)$ and $F_r: z\mapsto \Theta^+(r,z) + m_r(z)\Theta^-(r,z)$.
\end{Lem}
Recall that functions $\Theta^{\pm}$, $\Phi^{\pm}$ in Lemma \ref{l1} are defined in \eqref{eq32}.
\begin{Lem}\label{l2}
Let $\Hh$ be a singular nontrivial Hamiltonian on $\R_+$ whose spectral measure belongs to the Szeg\H{o} class $\sz$. Then $\K_{\Hh}$ is the non-negative non-increasing absolutely continuous function and we have 
$$
\lim_{r \to +\infty}\K_{\Hh}(r) = 0, \qquad \lim_{r \to +\infty}\K_{\widehat \Hh_r}(0) = \K_{\Hh}(0).
$$
\end{Lem} 
\begin{Lem}\label{l14}
Let $\Hh = \sth$ be a singular nontrivial Hamiltonian on $\R_+$ whose spectral measure belongs to the Szeg\H{o} class $\sz$. Then the functions $\I_{\Hh}$, $\Rr_{\Hh}$, $\K_{\Hh}$ are locally absolutely continuous and
\begin{align*}
\K'_{\Hh} &= 2\sqrt{h_1h_2 - h^2} - \I_\Hh h_1 - \frac{h_2}{\I_\Hh} + \frac{2\Rr_\Hh h}{\I_\Hh} - \frac{\Rr_\Hh^2 h_1}{\I_\Hh},\\
\frac{\I'_{\Hh}}{\I_{\Hh}} &= \I_{\Hh}h_1 - \frac{h_2}{\I_\Hh} + \frac{2\Rr_\Hh h}{\I_\Hh} - \frac{\Rr_\Hh^2 h_1}{\I_\Hh},  \\
\frac{\Rr'_{\Hh}}{\I_{\Hh}} &= 2 \Rr_{\Hh} h_1 - 2 h,
\end{align*} 
almost everywhere on $\R_+$. We also have $\K_{\Hh}(t) = \K_{\Hh^d}(t)$ for all $t \ge 0$ and the dual Hamiltonian $\Hh^d = J^* \Hh J$. 
\end{Lem} 

\medskip

\subsection{De Branges chains}\label{s23} Let $\Hh$ be a singular nontrivial Hamiltonian on $\R_+$. Define the Hilbert space
\begin{gather}\label{eq31}
L^2(\Hh) = \Bigl\{X \colon \R_+ \to \C^2 \colon \int_{\R_+}\bigr\langle\Hh(t)X(t), X(t)\bigr\rangle_{\C^2} \,dt < \infty\Bigr\}\Big/{\mathcal Ker}\,\Hh,\\
{\mathcal Ker}\,\Hh = \Bigl\{X\colon \;\Hh(t) X(t) = 0 \mbox{ for almost all } t \in \R_+\Bigl\}. \notag
\end{gather} 
As usual, the inner product in $L^2(\Hh)$ is defined by
$$
(X, Y)_{L^2(\Hh)} = \int_{\R_+}\langle\Hh(t)X(t), Y(t)\rangle_{\C^2}\,dt.
$$
An open interval $I \subset \R_+$ is called indivisible for $\Hh$ if there exists a vector $e \in \R^2$ such that $\Hh$ coincides with the rank-one operator $f \mapsto \langle f,e\rangle_{\C^2}e$ almost everywhere on $I$, and $I$ is the maximal open interval (with respect to inclusion) having this property. Let $\mathfrak{I}(\Hh)$ denote the set of all indivisible intervals of the Hamiltonian $\Hh$, and let $\Mi$ be the complement in $\R_+ \cup \{+\infty\}$ of the union of all intervals $I \in \mathfrak{I}(\Hh)$.  For $r \in \Mi$, define the subspace  
\begin{equation}\notag
H_r = \bigr\{X \in L^2(\Hh): \supp X \subset [0, r], \; X = x_I \mbox{ on } I\cap [0,r), \; I\in\mathfrak{I}(\Hh),\;x_I \in \C^2\Bigr\}.
\end{equation}
Let $\Theta$ be the first column of the matrix $M$ in \eqref{eq32}. Recall that $\langle \cdot , \cdot\rangle_{\C^2}$ denotes the inner product in $\C^2$,
$
\left\langle \left(\begin{smallmatrix}a_1\\ a_2\end{smallmatrix}\right) , \left(\begin{smallmatrix}b_1\\ b_2\end{smallmatrix}\right) \right\rangle_{\C^2} = a_1 \bar b_1 + a_2 \bar b_2.
$
Consider the mapping
\begin{equation}\label{wt}
\W_{\Hh}: X \to \frac{1}{\sqrt{\pi}}\int_{\R_+}\langle\Hh(t)X(t), \Theta(t, \bar z)\rangle_{\C^2}\,dt, \qquad z \in \C,  
\end{equation}
densely defined on $H_{+\infty}$ on functions with compact support. It is clear that $\W_{\Hh}X$ is the entire function with respect to $z$. The spectral measure of $\Hh$ could be characterized as the unique measure $\mu$ on $\R$ such that $(1 + x^2)^{-1} \in L^1(\mu)$ and the mapping $\W_{\Hh}$ is the isometric operator from $H_{+\infty}$ to $L^2(\mu)$. Uniqueness of the measure $\mu$ follows from the fact that the Hamiltonian is singular (relation \eqref{eq33} holds). More details can be found in Sections 8, 9 of \cite{Romanov} and in \cite{Winkler95}. 

\medskip

Given a nonzero measure $\mu$ on $\R$ with $(1+x^2)^{-1} \in L^1(\mu)$, we construct the de Branges chain associated to $\mu$ as follows. First, set $a = b =0$ and define $m$ by \eqref{hg}. Then use Krein--de-Branges theorem to construct a Hamiltonian $\Hh$ on $\R_+$ such that $m$ is the Weyl function for $\Hh$. Define subspaces $H_r \subset L^2(\Hh)$ and the transform $\W_{\Hh}$ in \eqref{wt}. The de Branges chain associated to $\mu$ consists of subspaces
$$
\B_r = \W_{\Hh}H_r, \qquad r \in \Mi.
$$ 
Note that the elements of $\B_r$ are entire functions. In fact, $\B_r$ is the Hilbert space equipped with the inner product inherited from $L^2(\mu)$. By construction, the subspaces $\B_r$ are isometrically embedded into $L^2(\mu)$. We also have isometric inclusions $\B_{r_1} \subset \B_{r_2}$ for $r_1 \le r_2$, $r_{1,2} \in \mathcal M$. Any regular de Branges space (a Hilbert space of entire functions having few natural properties) isometrically embedded into $L^2(\mu)$ coincides with one of spaces $\B_r$ in de Branges chain associated to $\mu$. This fact is a consequence of de Branges chain theorem, see Section 35 in \cite{dbbook}. In particular, it implies the uniqueness of de Branges chain of $\mu$.  

\medskip

In general, it is not known how to construct a nontrivial Hilbert space of entire functions isometrically embedded into $L^2(\mu)$ avoiding the usage of the inverse spectral theory for canonical Hamiltonian systems. However, for measures $\mu$ in the Szeg\H{o} class $\sz$ the situation simplifies. The following proposition is well-known for specialists.
\begin{Prop}\label{p1}
Let $\mu \in \sz$, and let $(\pw_{s}, \mu)$ be the completion with respect to the inner product of $L^2(\mu)$ of the linear space $\E_s$ of functions with smooth Fourier transform supported on $(-s, s)$. Then for every $s>0$ the Hilbert space $(\pw_{s}, \mu)$ belongs to the de Branges chain $\{\B_r\}$ associated to $\mu$. We have $(\pw_{s}, \mu) = \B_{r}$, where   
$r = \inf\{t: s = \int_{0}^{t}\sqrt{\det\Hh(\tau)}\,d\tau\}$. 
\end{Prop}

This result is a consequence of de Branges theorems 23, 35 in \cite{dbbook} and Krein formula for exponential type (Theorem 11 in \cite{Romanov}).

\medskip

Subspaces in de Branges chains admit a useful description in terms of Hermite-Biehler functions associated to canonical Hamiltonian systems. An entire function $E$ belongs to the Hermite-Biehler class if $|E(z)|>|E^\sharp(z)|$ for all $z \in \C^+$. Here and below we denote $f^\sharp(z) = \ov{f(\bar z)}$ for an entire function $f$. Let $\Hh$ be a singular nontrivial Hamiltonian on $\R_+$ and let $\mu$ be its spectral measure. We will say that a Hermite-Biehler function $E$ generates the subspace $\B_r$ of de Branges chain of $\mu$ if 
$$
\B_r = \left\{\mbox{entire } f: \; \frac{f}{E} \in H^2(\C^+), \; \frac{f^\sharp}{E} \in H^2(\C^+)\right\}, 
$$  
and, moreover, we have $\|f\|_{L^2(\mu)} = \|f\|_{L^2(|E|^{-2}\,dx)}$ for every $f \in \B_r$. As usual, $H^2(\C^+)$ denotes the standard Hardy space in the upper half-plane $\C^+$.  The standard choice for the  Hermite-Biehler function generating $\B_r$ is the function 
\begin{equation}\label{eq58}
E_r: z \mapsto \Theta^+(r,z) + i\Theta^-(r,z), \qquad z \in \C,
\end{equation} 
where $\Theta^{\pm}$ are the entries of the matrix $M$ in \eqref{eq32}, see Section 3 in \cite{Romanov}. If $E$ is any Hermite-Biehler function generating $\B_r$, then it is easy to see that the function
\begin{equation}\label{eq59}
k_{\B_r,\lambda}: z \mapsto -\frac{1}{2\pi i}\frac{E(z)\ov{E(\lambda)} - E^\sharp(z)\ov{E^\sharp(\lambda)}}{z - \bar \lambda}, \qquad z \in \C,
\end{equation} 
is the reproducing kernel of $\B_r$. This means that $k_{\B_r,\lambda} \in \B_r$ and we have
$f(\lambda) = (f, k_{\B_r,\lambda})_{L^2(\mu)}$ for every $f \in \B_r$. 

\medskip

Given a measure $\mu \in \sz$ and a number $s>0$, denote  by $(\pw_{[0, 2s]}, \mu)$ the set of functions $e^{isz} f$, where $f \in (\pw_{s}, \mu)$. Note that $(\pw_{[0, 2s]}, \mu)$ is the Hilbert space of entire functions with respect to the $L^2(\mu)$--inner product. Let $D_\mu$ be the Szeg\H{o} function \eqref{eq94} of $\mu$. Proposition \ref{p1} implies that $(\pw_{[0, 2s]}, \mu)$ can be identified with the subspace of the weighted Hardy space 
$$
H^2(w) = \left\{h = f D_\mu^{-1}, \;\; f \in H^2(\C^+)\right\}, \qquad w = |D_\mu|^2,
$$  
where the norm is defined by 
$$
\|h\|_{H^2(w)}^{2} = \int_{\R}|h(x)|^{2}w(x)\,dx.
$$
Choose any Hermite-Biehler function $E$ generating $\B_r = (\pw_s,\mu)$ and define the functions $P: z \mapsto e^{isz}E^\sharp$, $P^*: z \mapsto e^{isz}E$.
Using formula \eqref{eq59}, it is not difficult to see that
\begin{equation}\label{eq80}
k_{2s,\lambda}: z \mapsto -\frac{1}{2\pi i}\frac{P^{*}(z)\ov{P^{*}(\lambda)} - P(z)\ov{P(\lambda)}}{z - \bar \lambda}, \qquad z \in \C,
\end{equation}
is the reproducing kernel of $(\pw_{[0, 2s]}, \mu)$ at $\lambda \in \C$. By construction, we have $\mu = w\,dx + \mus$, hence $\|f\|_{H^2(w)} \le \|f\|_{L^2(\mu)}$ for every $f \in (\pw_{[0, 2s]}, \mu)$. It follows that $\|k_{2s,\lambda}\|_{L^2(\mu)}$ does not exceed the norm of the reproducing kernel of the space $H^2(w)$ at  $\lambda \in \C^+$. Since the reproducing kernel of $H^2(w)$ is given by 
$$k_{H^2(w),\lambda}: z \mapsto -\frac{1}{2\pi i}\frac{D_\mu^{-1}(z)\ov{D_\mu^{-1}(\lambda)}}{z - \bar \lambda},$$ 
this yields the useful inequality $|P^*(\lambda)|^2  - |P(\lambda)|^2 \le |D_\mu(\lambda)|^{-2}$ for every $\lambda \in \C^+$.

\medskip

\subsection{Reduction of Dirac systems to canonical Hamiltonian systems}\label{s24} It is well-known that Dirac systems could be rewritten as canonical Hamiltonian systems, see, e.g., Section 3 in \cite{Romanov}. For the reader convenience, we reproduce this observation. Consider the Dirac system
$$
JN'(t, z) + Q(t)N(t, z) = z N(t, z), \qquad  t \ge 0,\qquad N(0, z) = \idm,
$$ 
and let $N_0(t) = N(0, t)$, $t \ge 0$. Define $M = M(t,z)$ by $N(t,z) = N_0(t) M(t,z)$. Then we have
\begin{align*}
z N_0(t) M(t,z) 
&= J N'_0(t) M(t,z) + JN_0(t)J^* J M'(t,z) + Q N_0(t)M(t,z) \\
&= JN_0(t)J^* J M'(t,z) = (N^*_0(t))^{-1}JM'(t,z), 
\end{align*}
because $JAJ^* = (A^*)^{-1}$ for every real matrix $A$ with unit determinant. It follows that $M$ solves the Cauchy problem \eqref{cs} for the Hamiltonian $\Hh = N^*_0N_0$ on $\R_+$. In particular, for the potential $Q = \left(\begin{smallmatrix}0 & q\\ q  & 0\end{smallmatrix}\right)$ on $\R_+$ we have 
$$
N_0(t) = \begin{pmatrix}e^{-g(t)} & 0\\0 & e^{g(t)}\end{pmatrix}, \quad
\Hh(t) = \begin{pmatrix}e^{-2g(t)} & 0\\0 & e^{2g(t)}\end{pmatrix}, \quad t\ge 0,
$$
where $g(t) = \int_{0}^{t}q(s)\,ds$. Similarly, for the potential $Q = \left(\begin{smallmatrix}q & 0\\ 0  & -q\end{smallmatrix}\right)$ on $\R_+$ we have 
$$
N_0(t) = \begin{pmatrix}\cosh g(t) & \sinh g(t)\\\sinh g(t) & \cosh g(t)\end{pmatrix}, \quad
\Hh(t) = \begin{pmatrix}\cosh 2g(t) & \sinh 2g(t)\\\sinh 2g(t) & \cosh 2g(t)\end{pmatrix}, \quad t\ge 0.
$$
Next, let us show that the spectral measures of $\Hh$ and $\Di_Q$ coincide. Let $\Psi$ be the solution of \eqref{cp}, and let $\Theta$ be defined by \eqref{eq32}. From the relation $N = N_0 M$
we see that $\Psi(t, z) = N_0(t)\Theta(t, z)$ for all $t\ge 0$ and all $z \in \C$. Using the identity $\Hh = N_0^* N_0$, we get
$$
\int_{0}^{\infty}\langle\Hh(t) X(t), \Theta(t, \bar z)\rangle_{\C^2} = \int_{0}^{\infty}\langle N_0(t) X(t), \Psi(t, \bar z)\rangle_{\C^2}, 
$$
for every function $X \in L^2(\R_+, \C^2)$ with compact support. Moreover, the mapping $X \mapsto N_0 X$ is the unitary operator from $L^2(\Hh)$ to $L^2(\R_+, \C^2)$. It follows that operators $\F_Q: L^2(\R_+, \C^2) \to L^2(\mu)$,  $\W_\Hh: L^2(\R_+, \C^2) \to L^2(\mu)$ are isometric or not simultaneously. In other words, the spectral measures of $\Hh$ and $\Di_Q$ coincide. 

\medskip

\section{Khrushchev formula and proof of Theorem \ref{t2}}\label{s3}
Let $\Hh$ be a singular nontrivial Hamiltonian, and let $\{m_r\}_{r \ge 0}$ be the family of Weyl functions constructed in Section \ref{s2}. Recall that $m_r(i) = i \I_\Hh(r) + \Rr_\Hh(r)$ and $\I_{\Hh}(r) >0$ for every $r>0$. Define the Schur family $\{f_r\}_{r\ge 0}$ associated to $\Hh$ by
\begin{equation}\label{eq20}
m_r(z) = i \I_\Hh(r) \frac{1+B_i(z) f_r(z)}{1-B_i(z) f_r(z)} + \Rr_\Hh(r), \qquad z\in \C^+, \quad r\ge 0,
\end{equation}
where $B_i: z \mapsto \frac{z-i}{z+i}$ is the Blaschke factor in $\C^+$. By construction, the functions $f_r$ are analytic on $\C^+$. We also have
$$
\Im m_r(z) = \I_\Hh(r) \Re \left(\frac{1+B_i(z)f_r(z)}{1-B_i(z)f_r(z)}\right) = \I_\Hh(r) \frac{1 - |B_i(z) f_r(z)|^2}{|1 - B_i(z)f_r(z)|^2}.
$$
Since $\I_{\Hh}(r) >0$ for every $r \ge 0$, we see that $|B_i(z) f_r(z)| < 1$ for $z \in \C^+$. By Schwarz lemma,  this implies $|f_r(z)| \le 1$ for $z \in \C^+$. Thus, the family $\{f_r\}_{r \ge 0}$ consists of contracting analytic functions on $\C^+$. Let us represent the M\"obius transform $\tau: w \mapsto \I_\Hh(r) w + \Rr_{\Hh}(r)$ in the form 
\begin{equation}\label{eq36}
\tau: w \mapsto \frac{\sqrt{\I_{\Hh}(r)} w + \Rr_{\Hh}(r)/\sqrt{\I_{\Hh}(r)}}{ 0 \cdot w + 1/\sqrt{\I_{\Hh}(r)}}.
\end{equation} 
Let $\Theta^{\pm}$, $\Phi^{\pm}$ be the entries of the solution $M$ of Cauchy problem \eqref{cs}, see \eqref{eq32}. For $r>0$, define the entire functions $\tilde\Theta_r^{\pm}$, $\tilde\Phi_r^{\pm}$ by
\begin{equation}\label{eq18}
\begin{pmatrix}\tilde\Phi_r^{-} & \tilde\Phi_r^{+} \\ \tilde\Theta_r^{-} & \tilde\Theta_r^{+}\end{pmatrix} =
\begin{pmatrix}\Phi^{-}(r,\cdot) & \Phi^{+}(r,\cdot)\\ \Theta^{-}(r,\cdot) & \Theta^{+}(r,\cdot)\end{pmatrix}
\begin{pmatrix}\sqrt{\I_{\Hh}(r)} & \Rr_{\Hh}(r)/\sqrt{\I_{\Hh}(r)}\\ 0 & 1/\sqrt{\I_{\Hh}(r)}\end{pmatrix}.  
\end{equation}
Using \eqref{eq36} and \eqref{eq18}, we see that
\begin{equation}\label{eq19}
\frac{\tilde\Phi^-(r,z)w+\tilde\Phi^+(r,z)}{\tilde\Theta^-(r,z)w + \tilde\Theta^+(r,z)} =
\frac{\Phi^-(r,z)\tau(w) + \Phi^+(r,z)}{\Theta^-(r,z)\tau(w) + \Theta^+(r,z)}, \qquad w \in \C. 
\end{equation}
This formula will be used in the proof of Lemma \ref{l3} below. Relation \eqref{eq18} can be rewritten in the form
\begin{equation}\label{eq42}
\begin{pmatrix}\tilde\Theta_r^{+} & \tilde\Phi_r^{+} \\ \tilde\Theta_r^{-} & \tilde\Phi_r^{-}\end{pmatrix} =
\begin{pmatrix}1/\sqrt{\I_{\Hh}(r)} & \Rr_{\Hh}(r)/\sqrt{\I_{\Hh}(r)}\\ 0 & \sqrt{\I_{\Hh}(r)}\end{pmatrix} 
\begin{pmatrix}\Theta^{+}(r,\cdot) & \Phi^{+}(r,\cdot)\\ \Theta^{-}(r,\cdot) & \Phi^{-}(r,\cdot)\end{pmatrix}.
\end{equation}
The functions $\tilde\Theta_r^{\pm}$, $\tilde\Phi_r^{\pm}$ take real values on the real line and 
\begin{equation}\label{eq22}
\tilde\Theta_r^{+}(z)\tilde\Phi_r^{-}(z) - \tilde\Theta_r^{-}(z)\tilde\Phi_r^{+}(z) = 1, \quad z \in \C, 
\end{equation} 
as can be easily seen from \eqref{eq42} by calculating determinants and using the fact that $\det M(r,z) = \det M(0,z) = 1$.  Next, define the entire functions 
\begin{align}
\Et_r: z &\mapsto \tilde\Theta_r^+(z) + i \tilde\Theta_r^-(z),  &\Et_r^\sharp: z \mapsto \tilde\Theta_r^+(z) - i \tilde\Theta_r^-(z)\label{eq60}\\
\Ft_r: z &\mapsto \tilde\Phi_r^+(z) + i \tilde\Phi_r^-(z),  &\Ft_r^\sharp: z \mapsto \tilde\Phi_r^+(z) - i \tilde\Phi_r^-(z). \notag 
\end{align}
Since the functions $\Theta^{\pm}_r$ take real values on the real line, we have $|\Et_r(x)| = |\Et_r^\sharp(x)|$ for $x \in \R$. The same is true for the functions $\Ft_r$, $\Ft_r^\sharp$. 
\begin{Lem}\label{l4}
The measure $|\Et_r(x)|^{-2}\,dx$ is the spectral measure for $\widehat\Hh_r$. The function $\Et_r$ belongs to the Hermite-Biehler class and generates the subspace $\B_r = \W_\Hh H_r$ of de Branges chain associated to $\Hh$. 
\end{Lem}
\beginpf Let us apply assertion $(b)$ of Lemma \ref{l1} to the Hamiltonian $\widehat\Hh_r$. By definition, the Hamiltonian $(\widehat\Hh_r)_r$ coincides with the constant matrix
$$
(\widehat\Hh_r)_r = \begin{pmatrix}
1/\I_{\Hh}(r) & \Rr_{\Hh}(r)/\I_{\Hh}(r)\\
\Rr_{\Hh}(r)/\I_{\Hh}(r) & (\I_{\Hh}^2(r) + \Rr_{\Hh}^2(r))/\I_{\Hh}(r)
\end{pmatrix}.
$$
The Weyl function of this Hamiltonian equals 
$$
m_{(\widehat\Hh_r)_r} = i\I_{\Hh}(r) + \Rr_\Hh(r), 
$$
hence the spectral measure of $(\widehat\Hh_r)_r$ is $\I_{\Hh}(r)\,dx$. 
Since $\widehat\Hh_r = \Hh$ on $[0,r]$, we have $\hat F_r(z) = \Theta^+(r,z) + m_{(\widehat\Hh_r)_r}(z) \Theta^{-}(r,z)$ for the function $\hat F_r$ from Lemma \ref{l1} for $\widehat \Hh_r$. Then assertion $(b)$ of Lemma \ref{l1} says that for almost all $x \in \R$ we have 
\begin{equation}\label{eq95}
\hat w_r(x) = \frac{\I_{\Hh}(r)}{|\Theta^+(r,x) + (i\I_{\Hh}(r) + \Rr_\Hh(r)) \Theta^{-}(r,x)|^2} = \frac{1}{|\Et_r(x)|^2},
\end{equation}  
where $\hat w_r$ is the density of the absolutely continuous part of the spectral measure $\hat\mu_r$ of $\widehat\Hh_r$. Note that for every $x \in \R$ we have $\lim_{\eps \to 0}|F_r(x + i\eps)| > 0$ . Indeed, this follows from the fact that $\Theta^{\pm}(r, \cdot)$ are real analytic functions whose zero sets do not intersect: $\det M(r,z) = 1$. From assertion $(a)$ of Lemma~\ref{l1} we see that $\lim_{\eps \to 0}\Im \hat m_r(x + i\eps)$ exists and finite for every $x \in \R$, hence the singular part of $\hat\mu_r$ is zero. Formula \eqref{eq95} now says that $|\Et_r(x)|^{-2}\,dx$ is the spectral measure for $\widehat\Hh_r$. A lengthy but elementary calculation shows that 
$$
|\Et_r(z)|^2 - |\Et_r^\sharp(z)|^2 = 4 \Im(\Theta^+(z) \ov{\Theta^{-}(z)}) = |E_r(z)|^2 - |E_r^\sharp(z)|^2, \qquad z \in \C,
$$ 
where $E_r$ is defined by \eqref{eq58}. Since $E_r$ belongs to the Hermite-Biehler class, the above formula shows that the same is true for the function  $\Et_r$. Moreover, this formula implies that 
$$
\Et_r(z) \ov{\Et_r(\lambda)} - \Et_r^\sharp(z)\ov{\Et_r^\sharp(\lambda)} = E_r(z) \ov{E_r(\lambda)} - E_r^\sharp(z)\ov{E_r^\sharp(\lambda)}, \qquad z,\lambda \in \C. 
$$
Hence the reproducing kernels generated by $E_r$ and $\Et_r$ coincide, see \eqref{eq59}. It follows that $\Et_r$ generates the subspace $\B_r$. \qed 

\medskip

Next lemma is the analogue of Theorem 2 (page 173) in Khrushchev paper \cite{KH01}. 
\begin{Lem}\label{l3}
For every $r \ge 0$ and for almost all $x \in \R$ we have 
\begin{equation}\label{eq7}
|\Et_r(x)|^2 w(x) = \frac{1 - |\tilde f_r(x)|^2}{|1 - \theta_r(x) \tilde f_r(x)|^2},
\end{equation} 
where $\tilde f_r = B_i f_r$, and  $\theta_r: z \mapsto \frac{\Et^\ast_r(z)}{\Et_r(z)}$ is the inner function in $\C^+$.
\end{Lem}
\beginpf Assertion $(a)$ of Lemma \ref{l1} says 
$$
m(z) = \frac{\Phi^+(r,z) + m_r(z)\Phi^-(r,z)}{\Theta^+(r,z) + m_r(z)\Theta^-(r,z)}, \qquad z \in \C^+, \quad r \ge 0.
$$
From \eqref{eq20} and \eqref{eq19} we get
$$
m = \frac{\tilde\Phi_r^+ + i\tilde\Phi_r^-\frac{1+\tilde f_r}{1-\tilde f_r}}{\tilde\Theta_r^+ + i\tilde\Theta_r^-\frac{1+\tilde f_r}{1-\tilde f_r}} = \frac{(\tilde\Phi_r^+ + i \tilde\Phi_r^-) + (-\tilde\Phi_r^+ + i \tilde\Phi_r^-)\tilde f_r}{(\tilde\Theta_r^+ + i \tilde\Theta_r^-) + (-\tilde\Theta_r^+ + i \tilde\Theta_r^-)\tilde f_r}= \frac{\Ft_r - \Ft_r^\sharp\tilde f_r}{\Et_r - \Et_r^\sharp\tilde f_r} 
$$
in the upper half-plane $\C^+$. 
Taking into account \eqref{eq22} and the fact that $\tilde\Theta^{\pm}$, $\tilde\Phi^{\pm}$ are real on $\R$, we conclude that
$\Im \Ft_r\ov{\Et_r} = \Im (\tilde\Phi_r^+ + i \tilde\Phi_r^-)(\tilde\Theta_r^+ - i \tilde \Theta_r^-) = 1$ on the real line $\R$.
Analogously, $\Im \Ft_r^\sharp(x)\ov{\Et_r^\sharp(x)} = -1$ on $\R$. Next, the bounded analytic function $\tilde f_r$ has a non-tangential limit at almost every point $x \in \R$, and, moreover,
$$
\Im(\Ft_r^\sharp(x)\ov{\Et_r(x)}\tilde f_r(x) + \Ft_r(x)\ov{\Et_r^\sharp(x)\tilde f_r(x)}) = 0, 
$$
because $\Et^\sharp_r = \ov{\Et_r}$, $\Ft^\sharp_r = \ov{\Ft_r}$ on $\R$. It follows that for almost all $x \in \R$ we have
$$
\Im m(x) = \Im \frac{\Ft_r(x)\ov{\Et_r(x)}  + \Ft_r^\sharp(x)\ov{\Et_r^\sharp(x)}|\tilde f_r(x)|^2}{|\Et_r(x) - \Et_r^\sharp(x)\tilde f_r(x)|^2} = \frac{1 - |\tilde f_r(x)|^2}{|\Et_r(x) - \Et_r^\sharp(x)\tilde f_r(x)|^2}.
$$
It remains to note that $w(x) = \Im m(x)$ almost everywhere on $\R$, hence
$$
|\Et_r(x)|^2 w(x) = \frac{|\Et_r(x)|^2(1 - |\tilde f_r(x)|^2)}{|\Et_r(z) - \Et_r^\sharp(z)\tilde f_r(z)|^2} = \frac{1 - |\tilde f_r(x)|^2}{|1 - \theta_r(x) \tilde f_r(x)|^2},
$$
as required. Since $\Et_r$ belongs to the Hermite-Biehler class, the function $\theta_r$ is inner in $\C^+$.  \qed

\medskip

We are ready to prove Theorem \ref{t2}. It is interesting to note that the proof below is very close to the proof of Theorem 2.5   in Khrushchev paper \cite{KH01}. It seems that some other results from \cite{KH01} related to Schur functions also have analogues for canonical Hamiltonian systems.

\medskip

\noindent{\bf Proof of Theorem \ref{t2}.} For $r >0$, define the function $\Et_r$ by \eqref{eq60}. By Lemma \ref{l4}, $\Et_r$ belongs to the Hermite-Biehler class and generates the subspace $\B_r = \W_\Hh H_r$. Moreover, $|\Et_r(x)|^{-2}\,dx$ is the spectral measure for the Hamiltonian $\widehat\Hh_r$. Let $\hat m_r$ denote the Weyl function \eqref{mf} of $\widehat \Hh_r$. Since $\widehat \Hh_r$ coincides with $\Hh$ on $[0, r)$, the standard argument based on nesting circle analysis
gives $\lim_{r \to \infty}\Im \hat m_r(i) = \Im m(i)$, see Lemma 4.1 in \cite{BD2017} for more details. The last relation can be rewritten in the form
\begin{equation}\label{eq1}
\lim_{r \to \infty}\int_{\R}\frac{1}{|\Et_r(x)|^2}\,dP(x) = \int_{\R}w(x)\,dP(x),
\end{equation}
where we denoted $dP(x) = \frac{1}{\pi}\frac{dx}{x^2+1}$. Next, by Lemma \ref{l2} we have 
$$
\lim_{r \to \infty} \K_{\widehat \Hh_r}(0) = \K_{\Hh}(0).
$$ 
Using this fact, the definition of $\K_{\Hh}$, and relation \eqref{eq1}, we see that
\begin{equation}\label{eq2}
\lim_{r \to \infty}\int_{\R}\log\frac{1}{|\Et_r(x)|^2}\,dP(x) = \int_{\R}\log w(x)\,dP(x).
\end{equation}  
By Lemma \ref{l3}, we have
\begin{equation}\label{eq96}
\int_{\R}\log(|\Et_r|^2 w)\,dP 
= \int_{\R}\log \frac{1 - |\tilde f_r|^2}{|1 - \theta_r \tilde f_r|^2}\,dP 
= \int_{\R}\log (1 - |\tilde f_r|^2)\,dP,
\end{equation}
where we used the fact that $\int_{\R}\log|1 - \theta_r \tilde f_r|\,dP = \log|1 - \theta_r(i) \tilde f_r(i)| = 0$
by the mean value theorem for harmonic functions. Relations \eqref{eq2}, \eqref{eq96} now give us 
\begin{equation}\label{eq4}
\lim_{r \to \infty} \int_{\R}\log (1 - |\tilde f_r|^2)\,dP =  0.
\end{equation}
Denote $\log^+x = \max(\log x, 0)$. We have $\int_{\R}|\log g|\,dP \le 2 \int_{\R} \log^+ g \,dP$ for every function $g \in L^1(P)$ with $\int_{\R} \log g \,dP = 0$. Hence \eqref{eq4} implies
\begin{equation}\label{eq5}
\limsup_{r \to \infty}\int_{\R}\left|\log\frac{1}{|\Et_r|^2} - \log w\right|dP \le 2 \limsup_{r \to \infty}\int_{\R}\log^+\frac{1}{|\Et_r|^2 w}\,dP.
\end{equation}
Applying Lemma \ref{l3} once more, we obtain
\begin{align*}
\int_{\R}\log^+\frac{1}{|\Et_r|^2 w}\,dP 
&= \int_{\R}\log^+ \frac{|1 - \theta_r\tilde f_r|^2}{1 - |\tilde f_r|^2}\,dP
\le \int_{\R}\log\frac{(1 + |\tilde f_r|)^2}{1 - |\tilde f_r|^2}\,dP \\
&\le \int_{\R}\log \frac{1}{1 - |\tilde f_r|^2}\,dP + 2\int_{\R}|\tilde f_r|\,dP.  
\end{align*}
Using the fact that  $x \le \log\frac{1}{1-x}$ for $0 \le x < 1$ and $P(\R) = 1$, we can estimate 
$$
\left(\int_{\R}|\tilde f_r|\,dP\right)^2 \le \int_{\R}|\tilde f_r|^2\,dP \le \int_{\R}\log\frac{1}{1 - |\tilde f_r|^2}\,dP,  
$$
which tends to zero by \eqref{eq4}. We now see that \eqref{eq5} implies 
\begin{equation}\label{eq39}
\lim_{r \to \infty}\int_{\R}\left|\log\frac{1}{|\Et_r(x)|^2} - \log w(x)\right|dP(x) = 0.
\end{equation}
This completes the proof of the theorem. \qed 

\medskip

\section{Regularized Krein orthogonal entire functions}\label{s5}
Consider a Hamiltonian $\Hh = \sth$ on $\R_+$ such that $\det \Hh(t) = 1$  for almost all $t \in \R_+$. Let $\mu$ be the spectral measure of $\Hh$. Assume that $\mu = w\,dx +\mus$ belongs to the Szeg\H{o} class $\sz$. To simplify notations, we set $I = \I_\Hh$, $R = \Rr_{\Hh}$, $K = \K_{\Hh}$ for the rest of this section. For $r>0$, define entire functions $\Pt_{r}$, $\Pt_{r}^{*}$ by
\begin{equation}\label{eq77}
\Pt_{2r}: z \mapsto e^{irz - i\gphi(r)}\Et_{r}^\sharp(z), \qquad
\Pt^*_{2r}: z \mapsto e^{irz + i\gphi(r)}\Et_{r}(z), \qquad z \in \C,
\end{equation}
where $\gphi: r \mapsto \int_{0}^{r}\frac{R'(t)}{2I(t)}\,dt$ and $\Et_r$, $\Et_r^\sharp$ are given by \eqref{eq60}. Since $R$, $I$ are locally absolutely continuous and $I$ is strictly positive, the function $\gphi$ is correctly defined on $\R_+$ and is locally absolutely continuous. Functions $\Pt_r$, $\Pt_r^*$ could be regarded as regularized versions of Krein's orthogonal entire functions $P_r$, $P_r^*$ introduced by Krein in \cite{Kr81}, see \cite{Den06} for their modern theory.

\medskip

\begin{Lem}\label{l7}
For every $r>0$ the function $\Pt^{*}_{2r}$ is outer in $\C^+$. 
\end{Lem}
\beginpf Fix $r>0$. The function $\Et_r$ belongs the Hermite-Biehler class and 
\begin{equation}\label{eq81}
\limsup_{|z| \to \infty}\frac{\log|\Et_r(z)|}{|z|} = \limsup_{y \to +\infty}\frac{\log|\Et_r(iy)|}{y} = r,
\end{equation} 
see Section 6 in \cite{Romanov}. It follows that either $\Et_r = e^{-irz}$ or there exists a point $\lambda \in \C$, $\Im \lambda<0$, such that $\Et_r(\lambda) = 0$. In the first case we have $\Pt^{*}_{2r} = 1$ and there is nothing to prove. In the second case consider the reproducing kernel $k_{2r,\bar\lambda}$ of the space $(\pw_{[0,2r]},\mu)$. Using formula \eqref{eq80}, we see that 
$$
k_{2r,\bar\lambda}(z) = -\frac{1}{2\pi i}\frac{\Pt_{2r}^{*}(z) \ov{\Pt_{2r}^{*}(\bar\lambda)}}{z - \lambda}, \qquad z \in \C, 
$$
belongs to the weighted Hardy space $H^2(w)$. Hence $k_{2r,\lambda}/D_\mu \in H^2(\C^+)$ and we can use Smirnov-Nevanlinna factorization for functions of bounded type in~$\C^+$ (see, e.g., Theorem 9 in \cite{dbbook}) to find an inner function $\theta$ and an outer function $h$ such that $\Pt_{2r}^{*} =  \theta h$ in $\C^+$. Since $\Pt^{*}_{2r}$ does not vanish on $\{z \in \C: \; \Im z \ge 0\}$, the function $\theta$ has the form $\theta = e^{isz}$ for some $s \ge 0$. Moreover, we have
$$
\limsup_{y \to +\infty}\frac{\log|\Et_r(iy)|}{y} = \limsup_{y \to +\infty}\frac{\log|e^{ry}\Pt_{2r}^{*}(iy)|}{y} = r + \limsup_{y \to +\infty}\frac{\log|\theta(iy)|}{y} = r-s,
$$
because $\log|h(iy)| = o(y)$ is the Poisson extension of a function from $L^1(P)$. Comparing this formula with \eqref{eq81}, we see that $s = 0$. Hence, the function $\Pt_r^*$ is outer in $\C^+$, as required.
\qed

\begin{Lem}\label{l6}
The functions $I'/I$, $R'/I$, $(R'/I)^2(I h_1)^{-1}$, $1-Ih_1$, $1-1/Ih_1$ belong to the space $L^1(\R_+) + L^2(\R_+)$. We also have $K' \in L^1(\R_+)$. 
\end{Lem}
\beginpf The fact that that $K' \in L^1(\R)$ is the direct consequence of Lemma \ref{l2}. Relations in Lemma \ref{l14} can be rewritten in the form 
\begin{align}
-K' &= \left(I h_1 + \frac{1}{I h_1} - 2\right) 
+ \frac{1}{4}\left(\frac{R'}{I}\right)^2\frac{1}{I h_1},\label{eq61}\\
\frac{I'}{I} &= \left(I h_1 - \frac{1}{Ih_1}\right) - \frac{1}{4}\left(\frac{R'}{I}\right)^2\frac{1}{I h_1}, \label{eq62}\\
\frac{R'}{I} &= 2R h_1 - 2 h,
\end{align} 
almost everywhere on $\R_+$. Since $K' \in L^1(\R_+)$, and the functions 
$$
g = I h_1 + \frac{1}{I h_1} -2, \qquad \left(\frac{R'}{I}\right)^2\frac{1}{I h_1},
$$ 
are nonnegative, they belong to $L^1(\R_+)$. Denote $S_{1} = \{t \in \R_+: 1/2 \le I h_1 \le 2\}$.  
The function $g$ is comparable to $I h_1 + \frac{1}{I h_1}$ on $\R_+ \setminus S_{1}$, while on $S_{1}$ we have 
$$
c_1\left|1 - \frac{1}{I h_1}\right|^2 + c_1|1 - I h_1|^2 \le g \le c_2\left|1 - \frac{1}{I h_1}\right|^2 + c_2|1 - I h_1|^2
$$ 
for some positive constants $c_1$, $c_2$. Hence $1 - I h_1$, $1 - 1/I h_1$ are in $L^1(\R_+) + L^2(\R_+)$. From \eqref{eq61}, \eqref{eq62} we also see that $I'/I \in L^1(\R_+) + L^2(\R_+)$. Formula \eqref{eq61} implies that $R'/I \in L^2(S_{1})$. Moreover, on $S_2 = \R_+\setminus S_{1}$ we have
$$
\left(\int_{S_{2}}\left|\frac{R'(t)}{I(t)}\right|\,dt\right)^2 \le \int_{S_{2}}\left(\frac{R'(t)}{I(t)}\right)^2\frac{1}{Ih_1(t)}\,dt \cdot \int_{S_{2}}I h_1(t)\,dt < +\infty,
$$ 
where we used the fact that $g \in L^1(\R_+)$ is comparable to $I h_1 + \frac{1}{I h_1}$ on $S_2$. Thus, we have $R'/I \in L^1(\R_+) + L^2(\R_+)$ and lemma follows. \qed

\medskip

\noindent Next lemma concerns an analogue of the relation $\frac{\partial}{\partial r} P_r^*(z) = -A(r)P_{r}(z)$ (see Theorem 4.9 in \cite{Den06}) for functions $\Pt_r$, $\Pt_r^*$.
\begin{Lem}\label{l5}
For every $z \in \C$, the function $r \mapsto \Pt^*_r(z)$ is locally absolutely continuous on $\R_+$. Moreover, we have
\begin{equation}\label{eq43}
\tfrac{\partial}{\partial r} \Pt_{2r}^{*}(z) 
= - \tfrac{1}{2}(z-i)e^{2i\gphi(r)}\left(\tfrac{R'(r)}{I(r)} +i \tfrac{I'(r)}{I(r)}\right)\Pt_{2r}(z) + \tfrac{i}{2}z K'(r)\Pt_{2r}^{*}(z),
\end{equation}
for almost all $r > 0$.
\end{Lem}
\beginpf Relation \eqref{eq42} can we written in the form $\tilde\Theta = G \Theta$ , where the matrix-function $G$ is given by
\begin{equation}\label{eq23}
G(r) = \begin{pmatrix}1/\sqrt{I(r)} & R(r)/\sqrt{I(r)}\\ 0 & \sqrt{I(r)}\end{pmatrix}, \qquad r\ge 0. 
\end{equation}
Differentiating, we get 
$$
J \tilde\Theta' = J G' \Theta + z J G J^* \Hh \Theta = J G' G^{-1} \tilde \Theta + z J G J^* \Hh G^{-1} \tilde \Theta. 
$$
For the matrix $G$ we have
\begin{align}
G' &= \frac{1}{2}\frac{I'}{I}\begin{pmatrix}-1/\sqrt{I} & -R/\sqrt{I}\\ 0 & \sqrt{I}\end{pmatrix} + 
\begin{pmatrix}0 & R'/\sqrt{I}\\ 0 & 0\end{pmatrix}, \notag\\
G^{-1} &= \begin{pmatrix}\sqrt{I} & -R/\sqrt{I}\\ 0 & 1/\sqrt{I}\end{pmatrix} = (JGJ^*)^*, \label{eq25}\\
JG' G^{-1} &=-\frac{1}{2}\frac{I'}{I}\begin{pmatrix}0 & 1\\ 1 & 0\end{pmatrix} + \begin{pmatrix}0 & 0\\ 0 & R'/I\end{pmatrix},\notag\\
J G J^* \Hh G^{-1} &= \begin{pmatrix}I h_1 & -Rh_1 + h\\ -Rh_1 + h & (R^2 h_1 -2Rh + h_2)/I\notag\end{pmatrix}. 
\end{align}
Substituting these expressions into the formula for $J \tilde\Theta'$ and using identity $R'/I = 2Rh_1-2h$ from Lemma \ref{l14}, we obtain
\begin{align*}
\tfrac{\partial}{\partial r}\tilde\Theta^{-} &= \frac{1}{2}(I'/I+ z R'/I)\tilde \Theta^{-} - z Ih_1\tilde\Theta^{+},\\
\tfrac{\partial}{\partial r}\tilde\Theta^{+} &= -\frac{1}{2}(I'/I+ z R'/I)\tilde \Theta^{+} + (R'/I + z(R^2 h_1 -2Rh + h_2)/I)\tilde \Theta^{-}.
\end{align*}
It follows that
\begin{align*}
\tfrac{\partial}{\partial r} \Pt_{2r}^{*} 
= &i(z + \gphi'(r)) \Pt_{2r}^{*} + e^{irz + iu(r)}(\tfrac{\partial}{\partial r}\tilde\Theta^{+} + i\tfrac{\partial}{\partial r}\tilde\Theta^{-}),\\
= &i(z + \tfrac{R'}{2I})\Pt_{2r}^{*} - \tfrac{1}{2}e^{2i\gphi(r)}(\tfrac{I'}{I}+ z \tfrac{R'}{I})\Pt_{2r} \\
  &- iz e^{irz + i\gphi(r)}Ih_1\tilde \Theta^+ + e^{irz+ i\gphi(r)}(\tfrac{R'}{I} + z\tfrac{R^2 h_1 - 2Rh + h_2}{I})\tilde \Theta^{-}.
\end{align*}
Again by Lemma \ref{l14}, we have 
$$
\frac{R^2 h_1 - 2Rh + h_2}{I} = \frac{1}{4}\left(\frac{R'}{I}\right)^2\frac{1}{Ih_1} + \frac{1}{Ih_1} = - K' - Ih_1 + 2, 
$$
see also \eqref{eq61}. Using this identity, we obtain
\begin{align*}
- iz Ih_1\tilde \Theta^+ + (\tfrac{R'}{I} + z\tfrac{R^2 h_1 - 2Rh + h_2}{I})\tilde \Theta^{-}
&= - izIh_1\tilde \Theta^+ + (\tfrac{R'}{I} - zK' - zIh_1 + 2z)\tilde \Theta^{-}\\
&= - izIh_1\Et_r^\sharp(z) + (\tfrac{R'}{I} - zK' + 2z) \tilde \Theta^{-}.
\end{align*} 
Then relation $e^{irz}\tilde \Theta^{-} = \tfrac{1}{2i}(e^{-iu}\Pt_{2r}^{*} - e^{iu}\Pt_{2r})$ implies
\begin{align*}
\tfrac{\partial}{\partial r} \Pt_{2r}^{*} 
= &i(z + \tfrac{R'}{2I})\Pt_{2r}^{*} - \tfrac{1}{2}e^{2i\gphi(r)}(\tfrac{I'}{I}+ z \tfrac{R'}{I})\Pt_{2r} \\
  &- ize^{2i\gphi(r)}Ih_1\Pt_{2r} + \tfrac{1}{2i}(\tfrac{R'}{I} - zK' + 2z) (\Pt_{2r}^{*} - e^{2iu}\Pt_{2r})\\
= &\tfrac{iz}{2}K'\Pt_{2r}^{*} - \tfrac{1}{2}e^{2i\gphi(r)}(\tfrac{I'}{I}+ (z-i) \tfrac{R'}{I} + iz(2Ih_1 + K' - 2))\Pt_{2r}.	
\end{align*}
From \eqref{eq61}, \eqref{eq62} we get $2Ih_1 + K' - 2 = I'/I$, and formula \eqref{eq43} follows. \qed

\medskip

\begin{Lem}\label{l9}
Let $D_\mu$ be the Szeg\H{o} function of $\mu$. Then $\lim_{r \to +\infty}\Pt_r^{*}(z) = D_{\mu}^{-1}(z)$ for every $z \in \C^+$. Moreover,  
$\lim\limits_{r \to +\infty}\Pt_r(z) = 0$ and $\int_{\R_+}|\Pt_r(z)|^2\,dr < + \infty$ for $z \in \C^+$.
\end{Lem}
\beginpf Fix $z \in \C^+$. By Lemma \ref{l6} and Lemma \ref{l5}, we have  
\begin{equation}\label{eq71}
\tfrac{\partial}{\partial r}\Pt_{r}^{*}(z) = f_1(r)\Pt_r(z) + f_2(r)\Pt_{r}^{*}(z), 
\end{equation}
for some functions $f_{1} \in L^1(\R_+) + L^2(\R_+)$, $f_2 \in L^1(\R_+)$ depending on $z$. It follows that 
\begin{align*}
\tfrac{\partial}{\partial r}|\Pt_r^*(z)|^2 
&= 2\Re\left(\ov{\Pt_r^*(z)}\tfrac{\partial}{\partial r} \Pt_r^*(z)\right),\\
&= 2\Re\left(f_1(r)\Pt_r(z)\ov{\Pt_r^*(z)} + f_2(r)|\Pt_r^*(z)|^2\right),\\
&= 2\Re\left(f_1(r)\Pt_r(z)\ov{\Pt_r^*(z)}\right) - 2\Re f_2(r)|\Pt_r^*(z)|^2.
\end{align*}
Since $|\Pt_{r}(z)|^2 = e^{-2r \Im z }|\Pt_{r}^{*}(\bar z)|^2$, the above formula for $\bar z$ implies
\begin{align*}
\tfrac{\partial}{\partial r}|\Pt_{r}(z)|^2 
=&-2\Im z |\Pt_{r}(z)|^2 + 2e^{-2r\Im z}\Re\left(g_1(r)\Pt_r(\bar z)\ov{\Pt_r^*(\bar z)}\right) \\
&+ 2\Re g_2(r) e^{-2r\Im z}|\Pt_r^*(\bar z)|^2,
\end{align*}
for some functions $g_1 \in L^1(\R_+) + L^2(\R_+)$, $g_2 \in L^1(\R_+)$. Since $e^{-ir\bar z} \Pt_r(\bar z) = \ov{\Pt_r^*(z)}$ and $e^{irz} \ov{\Pt_r^*(\ov{z})} = 
\Pt_r(z)$, we have
\begin{equation}
\tfrac{\partial}{\partial r}|\Pt_{r}(z)|^2 
=-2\Im z |\Pt_{r}(z)|^2 + 2\Re\left(g_1(r)\Pt_r(z)\ov{\Pt_r^*(z)}\right) + 2\Re g_2(r)|\Pt_r(z)|^2. 	 \label{eq47}
\end{equation}
It follows that
\begin{align}
\tfrac{\partial}{\partial r}|\Pt_r^*(z)|^2  - \tfrac{\partial}{\partial r}|\Pt_{r}(z)|^2 
\ge&2\Im z |\Pt_{r}(z)|^2 -h_1(r)|\Pt_r(z)\Pt_r^*(z)| \notag\\
& - h_2(r)\left(|\Pt_r(z)|^2 + |\Pt_r^*(z)|^2\right), \label{eq44}
\end{align}
for some positive functions $h_1 \in L^1(\R_+) + L^2(\R_+)$, $h_2 \in L^1(\R_+)$. As we have seen in Section \ref{s23}, the function
\begin{equation}\label{eq75}
k_{2r,\lambda}: z \mapsto -\frac{1}{2\pi i}\frac{\Pt_{r}^{*}(z)\ov{\Pt_{r}^{*}(\lambda)} - \Pt_{r}(z)\ov{\Pt_{r}(\lambda)}}{z - \bar \lambda}, \qquad z \in \C,
\end{equation}
is the reproducing kernel of the space $(\pw_{[0,r]}, \mu)$ at $\lambda \in \C$. Hence for every $z \in \C^+$ the absolutely continuous function $\psi: r \mapsto |\Pt_{r}^*(z)|^2  - |\Pt_{r}(z)|^2$ increases and does not exceed $|D_\mu^{-1}(z)|^2$, where $D_\mu$ is the Szeg\H{o} function \eqref{eq94} of $\mu$. In particular, we have  $\psi' \in L^1(\R_+)$. Since $|\Pt_{0}(z)|>0$ and $\Pt_{r}(z)$ is continuous in $r$, there exists a number $r_0 > 0$ such that $\min_{0 \le r \le r_0}|\Pt_{r}(z)|^2 > 0$. On the other hand, for $r > r_0$ we have  $0<\psi(r_0) \le \psi(r) \le |\Pt_{r}^*(z)|^2$. Hence there is a constant $c>0$ such that $|\Pt_{r}^*(z)|^2 \ge c$ for every $r \ge 0$. Denote $\tilde h_{1,2} = h_{1,2}/2\Im z$. By the Hermite-Biehler property of $\Et_r$, we have $|\Pt_{r}(z)/\Pt_{r}^*(z)| \le 1$ for all $r \ge 0$. From \eqref{eq44} we get
\begin{equation}\label{eq74}
\frac{|\Pt_{r}(z)|^2}{|\Pt_{r}^*(z)|^2} \le \frac{\psi'(r)}{2c\Im z} + \tilde h_1(r)\left|\frac{\Pt_{r}(z)}{\Pt_{r}^*(z)}\right| + 
2\tilde h_2(r).	
\end{equation}
Consider the set $S = \{r \in \R_+: \; |\tilde h_1(r)| \ge 1\}$. We have $\tilde h_1 \in L^1(S)$, $\tilde h_1 \in L^2(\R_+\setminus S)$. Recall also that $\tilde h_2 \in L^1(\R_+)$ and $\psi' \in L^1(\R_+)$. Formula \eqref{eq74}, Cauchy inequality, and the bound $|\Pt_{r}(z)/\Pt_{r}^*(z)| \le 1$ imply the existence of a new constant $c>0$ such that
\begin{align*}
\int_{0}^{t} \frac{|\Pt_{r}(z)|^2}{|\Pt_{r}^*(z)|^2}\,dr 
&\le c + \|\tilde h_1\|_{L^1(S)} 
+ \|\tilde h_1\|_{L^2(\R_+ \setminus S)}\left(\int_{0}^{t} \frac{|\Pt_{r}(z)|^2}{|\Pt_{r}^*(z)|^2}\,dr \right)^{\frac{1}{2}}\!\!,\\
&\le c + c \left(\int_{0}^{t} \frac{|\Pt_{r}(z)|^2}{|\Pt_{r}^*(z)|^2}\,dr \right)^{\frac{1}{2}}\!\!,
\end{align*}
for every $t > 0$. It follows that 
\begin{equation}\label{eq72}
\frac{\Pt_{r}(z)}{\Pt_{r}^*(z)} \in L^2(\R_+)  \cap L^\infty(\R_+).
\end{equation}
Next, formula \eqref{eq43} for $z = i$ takes the form $\tfrac{\partial}{\partial r}\log \Pt_{2r}^*(i) = - K'(r)/2$. We have $K'\le 0$ almost everywhere on $\R_+$ by Lemma \ref{l2}. Since $\Pt_0^{*}(i) = 1/\sqrt{I(0)}$ is positive by construction, we see that $\Pt_r^{*}(i) > 0$ for all $r > 0$. Then Lemma \ref{l7} and Theorem \ref{t2} imply
\begin{equation}\label{eq73}
\lim_{r \to +\infty}\Pt_r^{*}(z) = D_{\mu}^{-1}(z), \qquad z \in \C^+.
\end{equation} 
Indeed, the functions $\Pt_r^{*}$, $D_{\mu}^{-1}$ are outer in $\C^+$, take positive values at $z = i$, and $\log |\Pt_r^{*}|$  tends to $\log |D_{\mu}^{-1}|$ in $L^1(P)$ as $r \to +\infty$ by Theorem \ref{t2}. Formulas \eqref{eq72} and \eqref{eq73} give us $|\Pt_r(z)|^2 \in L^1(\R_+) \cap L^\infty(\R_+)$. Moreover, we have $\frac{\partial}{\partial r}|\Pt_r(z)|^2 \in L^1(\R_+)$ by \eqref{eq47}. It follows that $\lim_{r \to +\infty} \Pt_r(z) = 0$ for every $z \in \C^+$. \qed

\medskip 

Next lemma is similar to Corollary 5.10 in \cite{KH01}.
	
\begin{Lem}\label{l15}
We have $\lim_{r \to+\infty}\|\Pt_{r}/(x + i)\|_{L^2(\mus)} = 0$.
Moreover if a set $E$ is such that $|\R\setminus E| = 0$, $\mus(E) = 0$, then 
$\lim_{r \to +\infty}\|(\Pt_{r}^* - \chi_{E} D_\mu^{-1})/(x+i)\|_{L^2(\mu)} = 0$. 
\end{Lem} 
\beginpf Fix the set $E$ from the statement of the lemma. From Theorem \ref{t2} and Jensen inequality we get
\begin{align}
1 &= \lim_{r \to \infty}\exp\left(\int_{\R}\log(|\Pt_{r}(x)|^2 w(x))\,dP(x)\right) \le \limsup_{r \to \infty}\int_{\R}|\Pt_{r}(x)|^2 w(x)\,dP(x) \notag\\
&\le \limsup_{r \to \infty}\frac{1}{\pi}\int_{\R}\frac{\chi_E(x)|\Pt_{r}(x)|^2}{1+x^2}\,d\mu(x) \le \limsup_{r \to \infty}\frac{1}{\pi}\int_{\R}\frac{|\Pt_{r}(x)|^2}{1+x^2}\,d\mu(x). \label{eq46}  
\end{align}
Recall that the reproducing kernel of the space $(\pw_{[0,r]}, \mu)$ is given by \eqref{eq75}. By Lemma \ref{l9}, we have
\begin{align}
\limsup_{r \to \infty}\left\|\frac{\Pt^*_{r}}{x+i}\right\|_{L^2(\mu)}^{2}	
&=\limsup_{r \to \infty}\frac{1}{|\Pt_{r}^{*}(i)|^2}\left\|\frac{\Pt_{r}^{*}(z)\ov{\Pt_{r}^{*}(i)} - \Pt_{r}(z)\ov{\Pt_{r}(i)}}{z + i}\right\|_{L^2(\mu)}^{2}\notag\\ 
&= \pi\limsup_{r \to \infty}\frac{|\Pt_{r}^{*}(i)|^2 - |\Pt_{r}(i)|^2}{|\Pt_{r}^{*}(i)|^2} = \pi. \label{eq78}
\end{align}
Since $|\Pt^{*}_{r}| = |\Pt_r|$ on $\R$, from here we see that inequalities in \eqref{eq46} are, in fact, equalities. It follows that
\begin{equation}\label{eq51}
\lim_{r \to+\infty}\|\Pt_{r}/(x + i)\|_{L^2(\mus)} = 0, \qquad \lim_{r \to+\infty}\|\Pt_{r}D_{\mu}\|_{L^2(dP)}^{2} = 1.
\end{equation}
Next, by the mean value theorem for the harmonic function $\Re(\Pt_{r}^{*}D_{\mu})$, we have 
\begin{align*}
\frac{1}{\pi}\int_{\R}\left|\frac{\Pt_{r}^*(x)}{x+i} - \frac{D_\mu(x)^{-1}}{x+i}\right|^2|D(x)|^2 \, dx
&=\int_{\R}\left|\Pt_{r}^*(x) D_{\mu}(x) - 1\right|^2 \, dP(x), \\
&=\|\Pt_{r}D_{\mu}\|_{L^2(dP)}^{2} + 1 - 2\Re(\Pt_{r}^*(i)D_{\mu}(i)).
\end{align*}
By Lemma \ref{l9} and \eqref{eq51}, the right hand side of the above formula converges to zero as $r \to +\infty$.
Using the fact that $\mu = |D_\mu(x)|^2\,dx + \mus$ and the first relation in \eqref{eq51}, we complete the proof of the lemma. \qed

\section{Proof of Theorem \ref{t1}}\label{s4}
\noindent We will use the following simple lemma, see, e.g, \cite{Den02b}\footnote{It worth be mentioned that $D_0 = -\Di_0$ for the free Dirac operator $D_0$ used in \cite{Den02b}. In particular, we have $e^{it\Di_0} = e^{-itD_0}$ for all $t\in \R$.}. 
\begin{Lem}\label{l17}
For every continuous function $f$ with compact support on $\R_+$ and for all $t>0$ large enough we have
$$
e^{-it\Di_0} \left(\begin{smallmatrix}f\\ 0\end{smallmatrix}\right) = \tfrac{1}{2}\left(\!\begin{smallmatrix}f(|s-t|)\\ if(|s-t|)\end{smallmatrix}\!\right), \qquad  
e^{-it\Di_0} \left(\begin{smallmatrix}0\\ f\end{smallmatrix}\right) = -\tfrac{i}{2}\left(\!\begin{smallmatrix}\sgn(s-t)f(|s-t|)\\ i\sgn(s-t) f(|s-t|)\end{smallmatrix}\!\right),
$$
for all $s \in \R_+$.
\end{Lem}
\noindent{\bf Proof of Theorem \ref{t1}.} Let $q_{1}$, $q_{2}$ be real functions on $\R_+$ such that $q_{1,2} \in L^{1}_{loc}(\R_+)$. Consider the Dirac operator $\Di_Q$ with the potential $Q = \left(\begin{smallmatrix}q_1 & q_2\\ q_2 & -q_1\end{smallmatrix}\right)$. Assume that the spectral measure $\mu$ of $\Di_Q$ belongs to the Szeg\H{o} class $\sz$. Using construction from Section \ref{s24}, define the Hamiltonian $\Hh = N_0^* N_0$ on $\R_+$ such that $\det\Hh = 1$ and the spectral measure of $\Hh$ coincides with $\mu$. Let $I =\I_{\Hh}$, $R = \Rr_{\Hh}$, $K = \K_{\Hh}$ be the functions from Section \ref{s22}, and let $G$ be defined by \eqref{eq23}. Consider the densely defined multiplication operator on $L^2(\R_+, \C^2)$ taking a continuous compactly supported vector-function $X \in L^2(\R_+, \C^2)$ into $NX$, where $N = J N_{0} J^* G^*\Sigma_{\gphi}$,
$$
\Sigma_\gphi(r) = \rotg, \qquad \gphi(r) = \frac{\Rr_{\Hh}'(r)}{2 \I_{\Hh}(r)}, \qquad r\ge 0.
$$
With a slight abuse of notation, we will denote this operator by the same letter~$N$. Let $\Psi$ be the solution of \eqref{cp}, and let $\Theta$ be defined by \eqref{eq32}, $\Tilde\Theta = G\Theta$. Since $N_0$ is a real matrix with $\det N_0 = 1$, we have $N^{-1} = J N_0^* J^*$. Using relation $\Psi = N_0 \Theta$ from Section \ref{s24}, we obtain
\begin{equation*}
\left\langle N(r) e, \Psi(r, x)\right\rangle_{\C^2} 
= \left\langle e, \Sigma_{\gphi}^{*}(r) G(r) N_0^{-1}(r) \Psi(r, x)\right\rangle_{\C^2}
= \bigl\langle e, \Sigma_{\gphi}^{*}(r)\tilde\Theta(r, x)\bigr\rangle_{\C^2},
\end{equation*}  
for every vector $e \in \C^2$ and all $r \ge 0$, $x \in \R$. Define $\Pt_{2r}$, $\Pt_{2r}^{*}$ by \eqref{eq77}. Consider a smooth function $f$ on $\R_+$ supported on $(0, a)$, $a>0$. Lemma \ref{l17} and the spectral theorem for $\Di_Q$ yield
\begin{align*}
\F_{Q} e^{it\Di_Q } N e^{-it\Di_0} \left(\begin{smallmatrix}f\\ 0\end{smallmatrix}\right) 
&=\frac{1}{\sqrt{\pi}}\int_{\R_+}\left\langle e^{it\Di_Q} N e^{-it\Di_0} \left(\begin{smallmatrix}f\\ 0\end{smallmatrix}\right), \Psi(s,x)\right\rangle_{\C^2}ds,\\ 
&= \frac{e^{itx}}{2\sqrt{\pi}} \int_{\R_+}\left\langle N(s) \left(\!\!\begin{smallmatrix}f(|s-t|)\\ if(|s-t|)\end{smallmatrix}\!\right), \Psi(s,x)\right\rangle_{\C^2}ds,\\
&=\frac{e^{itx}}{2\sqrt{\pi}}\int_{\R_+}f(|s-t|)\left\langle \Sigma_{\gphi}(s)\left(\begin{smallmatrix}1\\ i\end{smallmatrix}\right), \tilde\Theta(s,x)\right\rangle_{\C^2}ds,\\
&=\frac{e^{itx}}{2\sqrt{\pi}}\int_{\R_+}f(|s-t|)e^{-ixs}\Pt_{2s}^*(x)\,ds, 
\end{align*}
for $t>0$ big enough and all $x \in \R$. Next, expression $e^{ixt}\int_{\R_+}f(|s-t|)e^{-ixs}\Pt_{2s}^*(x)\,ds$ for $t>a$ equals
\begin{multline*}
\int_{-a}^{a}f(|s|)e^{-ixs}\Pt_{2s + 2t}^*(x)\,ds 
= \int_{-a}^{a}\left[\int_{-a}^{s}f(|\tau|)e^{-ix\tau}\,d\tau\right]' \Pt_{2s + 2t}^*(x)\,ds\\
= \Pt_{2a + 2t}^*(x)\int_{-a}^{a}f(|\tau|)e^{-ix\tau}\,d\tau + \int_{-a}^{a} \tfrac{\partial}{\partial s}\Pt_{2s + 2t}^*(x)\int_{-a}^{s} f(|\tau|)e^{-ix\tau}\,d\tau \,ds.
\end{multline*}
Denote the second summand in the last formula by $h_{t}(x)$. Let us show that $\lim_{t \to +\infty}\|h_t\|_{L^2(\mu)} = 0$. Since $f$ is smooth, we have
$$
\sup_{s \in (-a,a)}|\hat f_c(x, s)| \le \frac{c}{(1 + |x|)^2}, \qquad 
\hat f_c(x, s) = \int_{-a}^{s} f(|\tau|)e^{-ix\tau}\,d\tau, 
$$
for all $x \in \R$ and a constant $c$ depending only on $f$. Recall that $|\Pt| = |\Pt^*|$ on $\R_+$ by construction. By Lemma \ref{l5}, there exists a function $g \in L^1(\R) + L^2(\R_+)$ such that $\left|\frac{\partial}{\partial r}\Pt^*_r(x)\right| \le (1+|x|)g(x)|\Pt_r(x)|$ for $x \in \R$.
Hence one can use the generalized Minkowski inequality to estimate 
\begin{align}
\|h_t\|_{L^2(\mu)}^{2}
\le &c\int_{\R}\left(\int_{-a}^{a}g(s+t)\frac{|\Pt_{s+t}^{*}(x)|}{1+|x|}\,ds\right)^2d\mu(x) \notag\\
\le &c\left(\int_{-a}^{a}g(s+t)\,ds\right)^2 \sup_{s \in [-a, a]}\int_{\R}\frac{|\Pt_{s+t}^{*}(x)|^{2}}{(1+|x|)^2}\,d\mu(x).\label{eq48}
\end{align}
Formula \eqref{eq78} gives $\limsup_{r \to \infty}\int_{\R} \frac{|\Pt_{r}^{*}(x)|^2}{(1+|x|)^2}\,d\mu < +\infty$, hence
$\lim_{t \to \infty}\|h_t\|_{L^2(\mu)} = 0$. By Lemma \ref{l15}, functions $\Pt_{2a + 2t}^*\int_{-a}^{a}f(|\tau|)e^{-ix\tau}\,d\tau$ tend to $ \chi_{E} D_{\mu}^{-1}\int_{-a}^{a}f(|\tau|)e^{-ix\tau}\,d\tau$ in $L^2(\mu)$ as $t \to +\infty$. Hence there exists the strong limit in $L^2(\mu)$  
$$
\lim_{r \to +\infty} \F_Q e^{it\Di_Q}N e^{-it\Di_0} \left(\begin{smallmatrix}f\\ 0\end{smallmatrix}\right)
= \frac{\chi_{E} D_\mu^{-1}}{\sqrt{\pi}} \int_{0}^{\infty} f(t) \cos xt \,dt = \chi_{E} D_\mu^{-1}\F_0\left(\begin{smallmatrix}f\\ 0\end{smallmatrix}\right).  
$$
Similar arguments give the formula 
$$
\lim_{r \to +\infty} \F_Q e^{it\Di_Q}N e^{-it\Di_0} \left(\begin{smallmatrix}0 \\ f \end{smallmatrix}\right) 
= -\frac{\chi_{E} D_\mu^{-1}}{\sqrt{\pi}} \int_{0}^{\infty} f(t) \sin xt \,dt
= \chi_{E} D_\mu^{-1}\F_0 \left(\begin{smallmatrix}0 \\ f \end{smallmatrix}\right), 
$$
for smooth compactly supported functions $f$. Our aim now is to replace $N$ by the operator $X \mapsto \Sigma_\phi X$ for some function $\phi$. Let us define $\phi$ on $\R_+$ so that $\tilde N(r) = \Sigma_{\phi(r)} N(r)$ is positive definite for every $r \ge 0$. We are going to show that  
\begin{equation}\label{eq79}
\lim_{t \to + \infty}\|\tilde N S_t X - S_t X\|_{L^2(\R_+, \C^2)} = 0, \qquad (S_t X)(s) = X(t+s),
\end{equation}
for every continuous vector-function $X$ with compact support. By construction, we have $\det\tilde N = 1$ on $\R_+$. For $s \ge 0$, let $e_1(s)$, $e_2(s)$ be the unit eigenvectors of the matrix $\tilde N(s)$ corresponding to the eigenvalues $\lambda(s)$, $1/\lambda(s)$, correspondingly. Then for every vector $e = c_1e_1(r) + c_2e_2(r)$ in $\C^2$ we have 
$$
\|N e - e\|^2_{\C^2} = |c_1|^2(\lambda(s) - 1)^2 + |c_2|^2(1/\lambda(s) - 1)^2 \le (\lambda(s) - 1/\lambda(s))^2 \|e\|_{\C^2}^{2},
$$
where we used twice the elementary inequality $(x-1)^2 \le (x - 1/x)^2$ for $x>0$. Note that
\begin{align*}
(\lambda(s) - 1/\lambda(s))^2 
&= \trace\tilde N^2(s) - 2 = \trace N^*(s) N(s) - 2, \\
&= \trace G(s) J \Hh(s) J^* G^*(s) - 2.
\end{align*}
A calculation gives
$$
G J\Hh J^* G^* = \begin{pmatrix}\frac{h_2 - 2Rh + R^2h_1}{I} & Rh_1 - h\\ Rh_1-h & Ih_1\end{pmatrix} 
= \begin{pmatrix}1 & 0\\ 0 & 1\end{pmatrix} + \begin{pmatrix}\frac{1}{Ih_1}-1 + \frac{\gphi^2}{Ih_1} & \gphi\\ \gphi & Ih_1-1\end{pmatrix}.
$$
From Lemma \ref{l6} we see that $g = \trace G J \Hh J^* G^* - 2$ is a positive function in $L^1(\R_+)$. Hence 
\eqref{eq79} holds:
$$
\limsup_{t \to +\infty}\|\tilde N S_t X - S_t X\|_{L^2(\R_+, \C^2)}^{2} \le \limsup_{t \to +\infty}\int_{0}^{\infty}g(s)\|X(t+s)\|_{\C^2}^{2}\,ds = 0.
$$
It follows that for every smooth function $X$ with compact support the strong $L^2(\mu)$-limit  
$\lim_{t \to +\infty} e^{it\Di_Q}\Sigma_\phi e^{-it\Di_0} X$ exists and equals $\F_Q^{-1} \chi_{E} D_\mu^{-1}\F_0 X$. Since the operators $e^{it\Di_Q}\Sigma_\phi e^{-it\Di_0}$ are unitary, this implies the existence of the limit 
$$
W^m_-(\Di_Q, \Di_0) = \lim_{t \to +\infty} e^{it\Di_Q}\Sigma_\phi e^{-it\Di_0} = \F_Q^{-1} \chi_{E} D_\mu^{-1}\F_0
$$ 
in the strong operator topology. Existence of the wave operator 
$$
W^m_+(\Di_Q, \Di_0) = \F_Q^{-1} \chi_{E} \ov{D_\mu^{-1}}\F_0
$$ 
can be proved in a similar way. In the case where $q_1 = 0$, the Hamiltonian $\Hh$ generated by $Q$ has the diagonal form, see Section \ref{s24}. In particular, we have $\Rr_\Hh = 0$ by Lemma 2.2 in \cite{BD2017}. Hence our construction gives $\gphi = \phi = 0$, $\Sigma_\phi = \idm$, $W^{m}_{\pm}(\Di_Q, \Di_0) = W_{\pm}(\Di_Q, \Di_0)$ in this case. The theorem follows. \qed

\bigskip

Corollary \ref{c2} follows immediately from Theorem \ref{t1} and the following proposition.

\begin{Prop}\label{p3}
Let $q \in L^{1}_{loc}(\R_+)$ be a real function on $\R_+$, and let $Q = \left(\begin{smallmatrix}0 & q\\ q  & 0\end{smallmatrix}\right)$. Then the spectral measure of $\Di_Q$ belongs to $\sz$ if and only if $q$ satisfies \eqref{eq57}.   
\end{Prop}
\beginpf Consider the Hamiltonian  
$$
\Hh(t) = \begin{pmatrix}e^{-2g(t)} & 0\\0 & e^{2g(t)}\end{pmatrix}, \qquad g(t) = \int_{0}^{t}q(s)\,ds,
$$
on $\R_+$. Condition \eqref{eq57} for $q$ is equivalent to the condition \eqref{eq34} for $\Hh$. By Theorem~\ref{t4}, the spectral measure of $\Hh$ belongs to $\sz$. It remains to use the fact that the spectral measures of $\Hh$ and $\Di_Q$ coincide, see Section \ref{s24}. \qed 

\section{Appendix}
Here we prove Lemmas \ref{l1}, \ref{l2}, and \ref{l14} following the ideas of \cite{BD2017}. 

\medskip

\noindent{\bf Proof of Lemma \ref{l1}.} Assertions $(a)$, $(b)$ of Lemma \ref{l1} are the formulas $(2.13)$, $(2.14)$ in \cite{BD2017}, correspondingly. \qed

\medskip

\noindent{\bf Proof of Lemma \ref{l2}.} A straightforward calculation shows that the Weyl function of the constant nontrivial Hamiltonian 
$$
\Hh = \begin{pmatrix}c_1 & c \\ c & c_2\end{pmatrix}
$$
equals $m = i c_{1}^{-1}\sqrt{c_1 c_2 - c^2} + c/c_1$. Now let $\Hh$ be an arbitrary singular nontrivial Hamiltonian and let $\widehat\Hh_r$ be defined by \eqref{cook59}. Then the Weyl function of $(\hat \Hh_r)_r$ is
\begin{equation}\label{eq90}
m_{(\hat \Hh_r)_r} = i \I_{\Hh}(r) + \Rr_{\Hh}(r).
\end{equation}
Hence, we have $\J_{\widehat\Hh_{r}}(r) = \log c_{1}^{-1}\sqrt{c_1 c_2 - c^2} = - \log c_1 = \log \I_\Hh(r)$.
Next, let $F_{r}$, $G_r$ and $\hat F_{r}$, $\hat G_{r}$ be the functions from Lemma \ref{l1} for the Hamiltonians $\Hh$ and $\hat \Hh_r$, correspondingly. Note that $F_{r}(i) = \hat F_{r}(i)$, $G_{r}(i) = \hat G_{r}(i)$ by construction and formula \eqref{eq90}. It follows from assertion $(a)$ of Lemma \ref{l1} that $\hat m_r(i) = m_0(i)$, that is, 
$$
\I_{\widehat\Hh_r}(0) = \I_{\Hh}(0), \qquad \Rr_{\widehat\Hh_r}(0) = \Rr_{\Hh}(0).
$$ 
As in the proof of Lemma 2.5 in \cite{BD2017}, we have 
\begin{equation}\label{eq91}
\J_{\Hh}(r) = \J_{\Hh}(0) - 2\xi_{\Hh}(r) + 2\log|F_r(i)|, 
\end{equation}
where $\xi_{\Hh}: r \mapsto \int_{0}^{r}\sqrt{\det\Hh(t)}\,dt$. Similarly, 
$$
\J_{\widehat\Hh_r}(r) = \J_{\widehat\Hh_r}(0) - 2\xi_{\widehat\Hh_r}(r) + 2\log|\hat F_r(i)|.
$$
Since $\xi_{\Hh}(r) = \xi_{\widehat\Hh_r}(r)$ and $\hat F_r(i) = F_r(i)$, we have
\begin{align*}
\K_{\Hh}(r) 
&= \log \I_{\Hh}(r) - \J_{\Hh}(r) = \J_{\widehat\Hh_r}(r) - \J_{\Hh}(r),\\
&= \J_{\widehat\Hh_r}(0) - \J_{\Hh}(0) = \J_{\widehat\Hh_r}(0) - \log \I_{\widehat\Hh_r}(0) + \log \I_{\Hh}(0) - \J_{\Hh}(0), \\
&=-\K_{\widehat\Hh_r}(0) + \K_{\Hh}(0).
\end{align*}
The last formula can be rewritten in the form $\K_{\Hh}(0) = \K_{\Hh}(r) + \K_{\widehat\Hh_r}(0)$. Since the functions $\K_{\Hh}$, $\K_{\widehat\Hh_r}$ are nonnegative, we see that $\K_{\Hh}$ is nonincreasing. This fact and the semi-continuity of logarithmic integrals implies $\lim_{r \to +\infty} \K_{\widehat\Hh_r}(0) = \K_{\Hh}(0)$, or, equivalently, $\lim_{r \to +\infty} \K_{\Hh}(r) = 0$, see details in Lemma 4.1 of \cite{BD2017}. \qed

\medskip

\noindent{\bf Proof of Lemma \ref{l14}.} Assume first that the Hamiltonian $\Hh = \sth$ is continuously differentiable on $\R_+$. As in the proof of Lemma 2.7 in \cite{BD2017}, formula \eqref{cp} yields
\begin{equation}\label{eq92}
\left.\left(\!\begin{smallmatrix}\Theta^+(r,i)'&\Phi^+(r,i)'\\ \Theta^-(r,i)'&\Phi^-(r,i)'\end{smallmatrix}\!\right)\right|_{r=0}  = \left(\begin{smallmatrix}ih(0)&ih_2(0)\\-ih_1(0)&-ih(0)\end{smallmatrix}\right). 	
\end{equation}
Then \eqref{eq91} and the initial condition $M(0,i) = \idm$ give
\begin{align*}
\J_{\Hh}'(0) 
&= -2\xi'(0) + 2 \left.\Re\left(\frac{\Theta^+(r,i)' +m'_r(i)\Theta^-(r,i) + m_r(i)\Theta^-(r,i)'}{\Theta^+(r,i) +m_r(i)\Theta^-(r,i)}\right)\right|_{r=0},\\
&= -2\xi'(0) + 2\Re(ih(0) - im_0(i)h_1(0)),\\
&= -2\xi'(0) + 2 h_1(0) \I_{\Hh}(0),
\end{align*}
where the derivatives are taken with respect to $r$. The same relation for the Hamiltonian $\Hh_r$ in place of $\Hh$ shows that  $\J'_{\Hh}(r) = -2\xi'_{\Hh}(r) + 2 h_1(r) \I_{\Hh}(r)$ for all $r > 0$.	Similarly, differentiating relation 
$m_0(i) = \frac{G_r(i)}{F_r(i)}$ from Lemma \ref{l1} at $r = 0$ and using \eqref{eq92}, we obtain 
$0 = ih_2(0) + m'_0(i) - i m_0(i)h(0) - m_0(i)(ih(0) - i m_0(i)h_1(0))$.
As before, this gives
\begin{equation}\label{eq93}
0 = ih_2(r) + m'_r(i) - i m_r(i)h(r) - m_r(i)(ih(r) - i m_r(i)h_1(r))
\end{equation}
for all $r > 0$. By construction, we have $m_r(i) = i\I_{\Hh}(r) + \Rr_\Hh(r)$. Taking imaginary and real parts in \eqref{eq93}, we obtain
\begin{align*}
0&=h_2 + \I'_\Hh - 2\Rr_\Hh h + \Rr_\Hh^2 h_1 - \I_{\Hh}^{2}h_1,\\ 
0&=\Rr'_{\Hh}+2\I_{\Hh}h - 2\I_\Hh\Rr_\Hh h_1,
\end{align*}   
correspondingly. This two relations together with the definition of the entropy function $\K_{\Hh} = \log\I_{\Hh} - \J_\Hh$ and formula $\J'_{\Hh} = -2\xi'_{\Hh} + 2 h_1 \I_{\Hh}$ imply the formulas for $\K'_\Hh$, $\I'_\Hh/\I_\Hh$, and $\Rr'_\Hh/\I_\Hh$ in the case where $\Hh$ is smooth. These formulas in the general case then follow as in the proof of Lemma 2.7 in \cite{BD2017}. It remains to show that the spectral measure $\mu_r^d = w_r^d(x)\,dx + \mu_{\mathbf{s},r}^{d}$ for the Hamiltonian $\Hh_r^d: t \mapsto \Hh^d(t+r)$ belongs to $\sz$ and $\K_{\Hh^d}(r) = \K_{\Hh}(r)$ for every $r \ge 0$. This fact for diagonal Hamiltonians $\Hh$ is the part of Lemma 2.5 of \cite{BD2017}. However, the proof of this part actually does not uses the diagonal structure of $\Hh$ and hence works in our situation as well. \qed 

\vspace{0.5cm}

\bibliographystyle{plain} 
\bibliography{bibfile}

\def\cprime{$'$} \def\cprime{$'$} \def\cprime{$'$}
\begin{thebibliography}{10}

\bibitem{BD2017}
R.~V. Bessonov and S.~A. Denisov.
\newblock A spectral {S}zeg{\H o} theorem on the real line.
\newblock {\em Preprint arXiv:1711.05671}, 2017.

\bibitem{ChK}
M.~Christ and A.~Kiselev.
\newblock Scattering and wave operators for one-dimensional {S}chr\"odinger
  operators with slowly decaying nonsmooth potentials.
\newblock {\em Geom. Funct. Anal.}, 12(6):1174--1234, 2002.

\bibitem{dbbook}
L.~de~Branges.
\newblock {\em Hilbert spaces of entire functions}.
\newblock Prentice-Hall, Inc., Englewood Cliffs, N.J., 1968.

\bibitem{Den02b}
S.~A. Denisov.
\newblock On the existence of wave operators for some {D}irac operators with
  square summable potential.
\newblock {\em Geom. Funct. Anal.}, 14(3):529--534, 2004.

\bibitem{Den06}
S.~A. Denisov.
\newblock Continuous analogs of polynomials orthogonal on the unit circle and
  {K}re\u\i n systems.
\newblock {\em IMRS Int. Math. Res. Surv.}, pages 1--148, 2006, Art. ID 54517.

\bibitem{HSW}
S.~Hassi, H.~De~Snoo, and H.~Winkler.
\newblock Boundary-value problems for two-dimensional canonical systems.
\newblock {\em Integral Equations Operator Theory}, 36(4):445--479, 2000.

\bibitem{KK68}
I.~S. Kac and M.~G. Krein.
\newblock On the spectral functions of the string.
\newblock {\em Supplement II to the Russian edition of F.V. Atkinson, Discrete
  and continuous boundary problems}, 1968.
\newblock English translation: Amer. Math. Soc. Transl., (2) 103 (1974),
  19–102.

\bibitem{KH01}
Sergei Khrushchev.
\newblock Schur's algorithm, orthogonal polynomials, and convergence of
  {W}all's continued fractions in {$L^2(T)$}.
\newblock {\em Journal of Approximation Theory}, 108(2):161--248, 2001.

\bibitem{Kr81}
M.~G. Krein.
\newblock Continuous analogues of propositions on polynomials orthogonal on the
  unit circle.
\newblock {\em Dokl. Akad. Nauk SSSR (N.S.)}, 105:637--640, 1955.

\bibitem{LSb}
B.~M. Levitan and I.~S. Sargsjan.
\newblock {\em Sturm-{L}iouville and {D}irac operators}, volume~59 of {\em
  Mathematics and its Applications (Soviet Series)}.
\newblock Kluwer Academic Publishers Group, Dordrecht, 1991.
\newblock Translated from the Russian.

\bibitem{RSbook3}
Michael Reed and Barry Simon.
\newblock {\em Methods of modern mathematical physics. {III}}.
\newblock Academic Press, New York-London, 1979.

\bibitem{Romanov}
R.~Romanov.
\newblock Canonical systems and de {B}ranges spaces.
\newblock {\em preprint arXiv:1408.6022}, 2014.

\bibitem{Tep05}
A.~Teplyaev.
\newblock A note on the theorems of {M}. {G}. {K}rein and {L}. {A}.
  {S}akhnovich on continuous analogs of orthogonal polynomials on the circle.
\newblock {\em J. Funct. Anal.}, 226(2):257--280, 2005.

\bibitem{Winkler95}
H.~Winkler.
\newblock The inverse spectral problem for canonical systems.
\newblock {\em Integral Equations Operator Theory}, 22(3):360--374, 1995.

\end{thebibliography}


\enddocument